\documentclass[12pt]{article}
\usepackage{epsfig}
\usepackage{amsfonts,amsmath}

\newtheorem {thm}{Theorem}[section]
\newtheorem {problem}{Problem}[section]
\newtheorem {conjecture}[problem]{Conjecture}
\newtheorem {prob}[problem]{Problem}
\newtheorem {conj}[problem]{Conjecture}

\newtheorem {theo}[thm]{Theorem}

\newcommand{\bea}{\begin{eqnarray}}
\newcommand{\ba}{\begin{array}}
\newcommand{\bean}{\begin{eqnarray*}}
\newcommand{\ea}{\end{array}}
\newcommand{\eea}{\end{eqnarray}}
\newcommand{\eean}{\end{eqnarray*}}
\newcommand{\be}{\begin{equation}}
\newcommand{\ee}{\end{equation}}

\def\C{{\mathbb{C}}}
\def \K{{\mathcal K}}
\def\L{{\mathcal{L}}}

\def\F{{\mathcal{F}}}
\def\G{{\mathcal{G}}}

\def\1{{\mathbf 1}}

\def\R{{\mathbb R}}

\def\|{\, |\! | \, }

\def\al{\alpha}

\def\phi{\varphi}

\def\and{\, \mbox{ and } \,}

\def\b2b{\overline{\beta^2}}

\newcommand{\diam}{{\rm diam}}
\newcommand{\vol}{{\rm vol}}
\newcommand{\conv}{{\rm conv\;}}

\renewcommand{\and}{\hbox{ {\rm and} }}

\renewcommand{\C}{{\mathcal{C}}}

\newcommand{\RR}{{\mathbb R}}

\def\CHI{\mathchoice%
{\raise2pt\hbox{$\chi$}}%
{\raise2pt\hbox{$\chi$}}%
{\raise1.3pt\hbox{$\scriptstyle\chi$}}%
{\raise0.8pt\hbox{$\scriptscriptstyle\chi$}}}
\def\smalloplus{\raise1pt\hbox{$\,\scriptstyle \oplus\;$}}

\def\eps{\epsilon}

\newcommand{\bbR}{\mathbb{R}}
\newcommand{\bbZ}{\mathbb{Z}}
\newcommand{\Par}[1]{\left(#1\right)}
\newcommand{\Set}[1]{\{#1\}}
\newcommand{\su}{\subseteq}

\title  {Helly-type Problems}
\author{Imre B\'ar\'any and Gil Kalai}

\begin{document}

\maketitle

\begin{abstract}


In this paper, we present a variety of problems in the interface between 
combinatorics and geometry around the theorems of Helly, Radon,
Carath\'eodory, and Tverberg. Through these problems we describe the
fascinating area of Helly-type theorems, and explain some of its main themes
and goals.
\end{abstract}


\section{Helly, Carath\'eodory, and Radon theorems}

In this paper, we present a variety of problems in the interface between 
combinatorics and geometry around the theorems of Helly, Radon,
Carath\'eodory, and Tverberg.


Helly's theorem~\cite{Helly:1923wr} asserts that for
a family $\{K_1,K_2,\ldots, K_n\}$ of convex sets in $\mathbb R^d$ where $n \ge d+1$,
if every $d+1$ of the sets have a point in common, then all of the sets have a point in common.
The closely related Carath\'eodory theorem~\cite{Carath1907} states that for $S \subset \mathbb R^d$,
if $x \in \conv S$, then $x \in \conv R$ for some $R \subset S$, $|R| \le d+1$.

The more general colorful Carath\'eodory theorem~\cite{Barany:1982va} says the following. 
Let $S_1$, $S_2, \ldots ,S_{d+1}$ be $d+1$ sets (or colors if you wish) in $\mathbb R^d$.
Suppose that $x \in \bigcap_{i=1}^{d+1}\conv S_i$. Then there is a {\sl transversal}  $T=\{x_1,\ldots,x_{d+1}\}$
of the system $S_1,\ldots,S_{d+1}$, meaning
that $x_1 \in S_1, x_2 \in S_2, \ldots, x_{d+1} \in S_{d+1}$ such that $x \in \conv T$.
A transversal is also called a {\sl rainbow set} when $S_1,\ldots,S_{d+1}$ are considered as colors.
The uncolored version, that is, when $S_1=S_2=\ldots =S_{d+1}$, is the classic result of Carath\'eodory.
There is a closely related colorful version of Helly's theorem due to Lov\'asz that appeared in~\cite{Barany:1982va}.

Tverberg's theorem~\cite{Tverberg:1966tb} states the following: Let $x_1,x_2,\ldots, x_m$ be points 
in $\mathbb R ^d$ with $m \ge (r-1)(d+1)+1$. Then there is a partition  $S_1,S_2,\ldots, S_r$ of $\{1,2,\ldots,m\}$ 
such that $\bigcap _{j=1}^r\conv \{x_i: i \in S_j\} \ne \emptyset$. This was a conjecture
by Birch who also proved the planar case in a slightly different form. The bound of $(r-1)(d+1)+1$ in the theorem is 
sharp as can easily be seen from the configuration of points in a sufficiently general position.

The case $r=2$ is Radon's theorem~\cite{Radon:1921vh}, another classic from 1921,
which was used by Radon to prove Helly's theorem. Helly's original proof (published later) was based on
a separation argument.
Sarkaria~\cite{Sarkaria:1992vt} gave a simple proof of Tverberg's
theorem based on the colorful Carath\'eodory theorem.

This paper describes the fascinating area of Helly-type theorems, and explains some of its main themes
and goals through a large and colorful bouquet of problems and conjectures.
Some of these problems are very precise 
and clear-cut, for instance, Sierksma's conjecture (Conjecture~\ref{sierksma}), the cascade conjecture (Conjecture~\ref{cascade}), and Problem~\ref{p:volum} about volumes of intersections. 
Some of them are rather vague, for instance, Problem~\ref{d-repr} about intersection patterns of Euclidean convex sets, and Problem \ref{mutual} about the mutual position of convex sets, and Problem \ref{vague} about topological conditions for the existence of Tverberg partitions. 
We hope to see the answers to many of the questions presented here in the near future.
Often, results from convexity give a simple and strong
manifestation of theorems from topology: Helly's theorem manifests the nerve theorem from algebraic topology,
and Radon's theorem can be regarded as an early ``linear'' version of the Borsuk--Ulam theorem.
One of our main themes is to further explore these connections to topology.
Helly-type theorems also offer
complex and profound combinatorial connections and applications that represent a
second theme of this paper.

For a wider perspective and many other problems we refer the reader to survey papers
by Danzer--Gr\"unbaum--Klee \cite {Danzer:1963ug}, Eckhoff \cite {Eckhoff:1979bi}
and \cite {Eck93survey}, Tancer \cite {Tancer:2013iz}, De Loera--Goaoc--Meunier--Mustafa \cite{de2017discrete},
and the forthcoming book of B\'ar\'any \cite{Bar2021}.


Here is a quick summary of the paper. Section \ref {s:h1}
defines the nerve that records the intersection pattern of
convex sets in $\mathbb R^d$, describes some of its combinatorial and topological properties, and considers
various extensions of Helly's theorem, such as the fractional Helly theorem, which asserts
that if a fraction $\alpha$ of all sets in a family of convex sets have a 
non-empty intersection, then there is
a point that belongs to a fraction $\beta (\alpha,d)$ of the sets in the family. 
Section \ref {s:h2} considers various refinements and generalizations of Helly theorems such as the
study of dimensions of intersections of convex sets, and the study of Helly-type theorems for unions of convex sets.
Section \ref{s:t1} presents various extensions and refinements of Tverberg's theorem, starting with Sierksma's
conjecture on the number of Tverberg partitions. Section \ref{s:t2} studies the cascade conjecture about the dimensions of the Tverberg points and considers
several connections with graph theory including a speculative connection with the four-color theorem.
Section \ref{s:t3} deals with other Tverberg-type problems. Section \ref{s:c} brings
problems related to Carath\'eodory theorem and weak-epsilon nets, and Section
\ref {s:ct} gives a glance at common transversals; rather than piercing a family of sets by
a single point or a few points we want to stab them with a single or a few $j$-dimensional affine spaces.
Final conclusions are drawn in the last section.

\section {Around Helly's theorem} 
\label {s:h1}

\subsection {Nerves, representability, and collapsibility}

We start this section with the following basic definition:
for a finite collection of sets $\F = \{K_1,K_2,\dots,K_n\}$,
the {\sl nerve} of $\F$ is a simplicial complex defined by
$${\cal N}({\F})
=\{S \subset [n]: \bigcap _{i \in S}K_i \ne \emptyset\}.$$
Helly's theorem can be seen as a statement about nerves of
convex sets in $\mathbb R^d$, and nerves come to play in many extensions and refinements of Helly's theorem.

A {\sl missing face} $S$ of a simplicial complex $\K$ is a set of vertices of $\K$ that is not
a face but every proper subset of $S$ is a face. Helly's theorem asserts that a $d$-{\sl representable} complex does
not have a missing face with more than $d+1$ vertices.

A simplicial complex is $d$-{\sl representable} if it is the
nerve of a family of convex sets in $\R^d$.

\begin {prob}\label{d-repr}
Explore $d$-representable simplicial complexes.
\end {prob}

We refer the reader to the survey on $d$-representable complexes by Tancer \cite {Tancer:2013iz}.

Let $\K$ be a simplicial complex. A face $F \in \K$ is {\sl free} if it is contained in a unique maximal face.
An elementary $d$-{\sl collapse} step is the removal from $\K$ of a free face $F$ with at most $d$
vertices and all faces containing $F$. A simplicial complex is $d$-{\sl collapsible}
if it can be reduced to the empty complex by a sequence of elementary $d$-collapse steps.
Wegner proved \cite{Wegner:1975eo}
that every $d$-representable complex is $d$-collapsible. The converse does not hold even for $d=1$:
1-representable complexes are
the clique complex of interval graphs and 1-collapsible complexes are the clique complexes of chordal graphs.
 
Here, a clique complex of a graph $G$ is a simplicial complex whose faces correspond
to the sets of vertices of complete subgraphs of $G$. Chordal graphs are graphs
with no induced cycles of length greater than three.
Intersection patterns of intervals (which are the same as 1-representable complexes)
were completely characterized by Lekkerkerker and Boland \cite{LeBo:1962}. They proved that interval graphs 
are characterized by being chordal graphs with the additional property that among every three vertices,
one is a vertex or adjacent to a vertex in any path between the other two. They also described interval
graphs in terms of a list of forbidden induced subgraphs.


\subsection {The upper bound theorem}

For a finite collection of 
sets ${\cal F} = \{K_1,K_2,\dots,K_n\}$, $n \ge d+1$,
in $\mathbb R^d$, let ${\cal N}={\cal N}({\cal F})$
be the nerve of $\cal F$. We put $f_k({\cal N})=|\{S \in {\cal N}: |S|=k+1\}|$.
(The vector $(f_0({\cal N}), f_1({\cal N}),\dots)$ is called the $f$-{\sl vector} of ${\cal N}$, and is sometimes referred to also
as the $f$-vector of $\cal F$ and written as $f(\F)$.)
Helly's theorem states that if $f_{n-1}({\cal N})=0$ then $f_d({\cal N}) < {{n} \choose {d+1}}$, or , with the $f(\F)$ notation, $f_{n-1}(\F)=0$ implies $f_d(\F) < {{n} \choose {d+1}}$

A far-reaching extension of Helly's theorem was conjectured by Katchalski and Perles and
proved by Kalai \cite{Kalai:1984bg} and Eckhoff \cite {Eckhoff:1985hr}.

\begin{theo}[upper bound theorem]\label{t:upper_bound}
 Let $\F$ be a family of $n$ convex sets in $\R^d$, and suppose that every $d+r+1$ members of $\F$ have an empty intersection.
 Then, for $k=d,\ldots, d+r-1,$
 \[
  f_k({\cal N}(\F)) \leq \sum_{j=k}^{d+r-1}\binom{j-d}{k-d}\binom{n-j+d-1}{d}.
 \]
\end{theo}

The theorem provides best upper bounds for $f_d(\F),\ldots,f_{d+r-1}(\F)$ in terms of $f_0(\F)$ provided $f_{d+r}(\F)=0.$
The proofs rely on $d$-collapsibility. 
There is a simple case of equality: the family consists of $r$ copies of $\mathbb R^d$ and $n-r$
hyperplanes in general position. Theorem \ref {t:upper_bound} is
closely related to the upper bound theorem for convex polytopes of Peter McMullen~\cite{McM70}. 
In fact, a common proof was given by Alon and Kalai in~\cite{Alon:1995fs}.

\begin {prob}
\label {p:ubte} Study cases of equality for the upper bound theorem.
\end {prob}

A place to start would be to understand $2$-representable complexes
$\K$ with $f_3(\K)=0$ and $f_2(\K)={{n-1} \choose {2}}$. 

Theorem \ref {t:upper_bound} implies the sharp version of the fractional Helly theorem 
of Katchalski and Liu~\cite {Katchalski:1979vt}. The sharp version is due to Kalai~\cite{Kalai:1984bg}.

\begin {theo}
  \label {t:fht}
  Let $\K$ be a $d$-representable complex. If $f_d(\K) \ge \alpha {{n \choose {d+1}}}$, then $\dim (\K) \ge \beta n$,
  where $\beta = \beta (d,\alpha)=1- (1-\alpha)^{1/d+1}$. 
\end {theo}

In other words, if $\F=\{K_1,\ldots,K_n\}$ is a family of convex sets in $\R^d$ ($n\ge d+1$)
and at least $\alpha{n \choose d+1}$ of the $d+1$ tuples in $\F$ intersect, 
then $\F$ contains an intersecting subfamily of size $\beta n$. This is a result of central
importance around Helly's theorem. The existence of $\beta (\alpha)$ is
referred to as the {\sl fractional Helly property}, and if $\beta \to 1$ when $\alpha \to 1$
this is referred to as the {\sl strong} fractional Helly property.

We note that a complete characterization of $f$-vectors of $d$-representable complexes was
conjectured by Eckhoff and proved by Kalai \cite{Kalai:1984isa,Kalai:1986hoa}.

\subsection {Helly numbers and Helly orders}

\label{s:hnho}

It is useful to consider the following abstract notions of Helly numbers and Helly orders.
Let ${\cal F}$ be a family of sets.
The {\it
Helly number} $~{\rm h} ({\cal F})$ of ${\cal F}$ is the minimal positive
integer $h$ such that if a finite subfamily ${\cal K} \subset {\cal F}$
satisfies $\bigcap {\cal K}' \neq \emptyset$ for all ${\cal K}'
\subset {\cal K}$ of cardinality $\leq h$, then $\bigcap {\cal K}
\neq \emptyset$.
The {\it Helly order} $~{\rm h_o} ({\cal F})$ of ${\cal F}$ is the minimal positive
integer $h$ such that if a finite subfamily ${\cal K} \subset {\cal F}$
satisfies

(1) every finite intersection of sets in ${\cal K}$ belongs to ${\cal F}$
and

(2) $\bigcap {\cal K}' \neq \emptyset$ for all ${\cal K}'
\subset {\cal K}$ of cardinality $\leq h$,

\smallskip
\noindent then $\bigcap {\cal K} \neq \emptyset$. Of course, when we consider
families of sets closed under intersection, the Helly number
and the Helly order coincide.
So, for example, the topological Helly theorem, to be mentioned next in Section \ref {s:th},
asserts that the Helly order of topologically trivial sets in $\mathbb R^d$
is $d+1$, and Amenta's theorem (Section \ref{s:h-u}
below) asserts that
the family of unions of $k$ pairwise disjoint convex sets
in $\mathbb R^d$ has the Helly order $k(d+1)$.

Let ${\cal F}$ be a family of sets.
The {\it fractional Helly number} $~{\rm g} ({\cal F})$ of ${\cal F}$ is the minimal positive
integer $g$ such that there is a function $f(\alpha)>0$, defined for $\alpha >0$, with the following property: 
for every family ${\cal K} \subset {\cal F}$ of cardinality $n$,
if at least $\alpha{n \choose g}$ of the $g$-tuples in $\F$ intersect, 
then $\F$ contains an intersecting subfamily of size $f(\al) n$.  


\subsection {Topological Helly theorem and Leray complexes}
\label {s:th}

Helly himself proved a topological version of his theorem~\cite{Helly:1930hk}.
A {\sl good cover} is a family of compact subsets of $\mathbb R^d$
such that every intersection of sets in the family is either empty or
topologically trivial. (By ``topologically trivial''
we mean ``contractible,'' but it is sufficient to assume that all homology groups vanish.)

\begin {theo}[topological Helly]
  If in a good cover of $n$ subsets of $\mathbb R^d$, $n \ge d+1$,
  every intersection of $d+1$ sets is non-empty, then the intersection of
  all the sets in the family is non-empty. 
\end {theo}

A simplicial complex $\K$ is $d$-{\sl Leray} if $H_i(\K')=0$ for every
induced subcomplex $\K'$ of $\K$ and for every $i \ge d$.
The well-known nerve theorem from algebraic topology asserts that if $\cal K$ is a
finite family of sets that form a good cover then the
nerve of $K$ is topologically equivalent to $\bigcup {\cal K}$. (The notion of ``topologically equivalent''
corresponds to the notion of ``topologically trivial'' in the definition of good covers.)
It follows from the homological version of the nerve theorem that $d$-representable complexes are $d$-Leray.
It is also easy to see that $d$-collapsible complexes are $d$-Leray.

\smallskip
{\bf Remark:} The nerve theorem played an important role in algebraic topology in the '40s and '50s, e.g., in showing
that the de Rham homology coincides with other notions of homology.
Helly's topological theorem is remarkable since it came earlier than these
developments. In Section \ref {s:toptve} we will mention that
Radon's theorem can be seen as an early incarnation of the Borsuk--Ulam theorem in topology. 
Topological extensions of Helly-type theorems are an important part of the theory.
Often, such extensions are considerably more difficult to prove, but in a few cases the topological proofs
are the only known ones even for the geometric results. There are also a few cases where natural
topological extensions turned out to be incorrect. The survey paper~\cite{de2017discrete} of De Loera, Goaoc, Meunier, and Mustafa 
emphasizes connections with combinatorial theorems closely related to the Brouwer fixed-point theorem, starting with
the Sperner lemma and the Knaster--Kuratowski--Mazurkiewicz theorem.

\smallskip
A general problem is the following.

\begin {prob}
\label {meta-problem}
(i) Find finer and finer topological and combinatorial properties of $d$-representable complexes.

(ii) Extend Helly-type theorems to good covers, Leray complexes, and beyond.

(iii) Find weaker topological conditions that suffice for the topological Helly theorem to hold. 
  
\end {prob}

There is much to say about part (ii) of Problem \ref {meta-problem}.
In several cases the way to go about it is to extend
properties of $d$-representable complexes to $d$-Leray complexes. We will come back
to such extensions later but we note that the upper bound theorem
(Theorem \ref{t:upper_bound}), as well as the full characterization of their $f$-vectors,
extends to $d$-Leray complexes; ~\cite{kalai:2002ij}. 
This is also closely related to Stanley's characterization~\cite{Stanley:1975} of $f$-vectors
of Cohen--Macaulay complexes.

\begin{figure}[h!]
\centering
\includegraphics[scale=0.8]{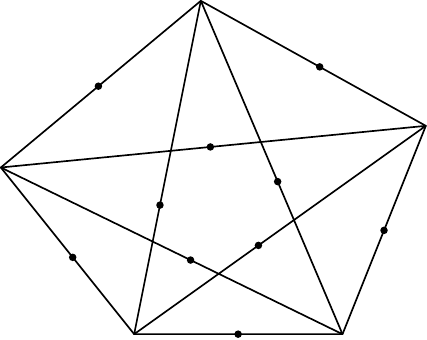}
\caption{$TK_5$}
\label{fig:Tk5}
\end{figure}

Regarding part (i) of Problem \ref {meta-problem}, we first note a very easy connection
with embeddability:
if $G$ is a graph we denote by $TG$ the graph where for every edge $e$,
we add a new vertex $v_e$  that is adjacent
to the endpoints of $e$, and remove $e$ itself; see Figure~\ref{fig:Tk5}.
If $G$ is not planar, e.g., when $G=T_5$,
then $TG$ is not 2-representable. (Note however that $K_n$ itself is 2-representable
for every $n$.)

White~\cite{White2020} 
defined the class of {\it $d$-Matou{\v s}ek} simplicial complexes that are  related to topological invariants for embeddability  
as follows. Let $\K$ be an abstract simplicial complex with vertices $V(\K)=[n]=\Set{1,2,\ldots,n}$.
We define the dual simplicial complex $K'$, with
vertices $V(\K')=\Set{J\in \K\ :\ J\text{ is inclusion maximal}}$,
and faces $\K'=\Set{\alpha\su V(\K')\ :\ \bigcap_{J\in\alpha}J\neq\emptyset}$.\\
We say that $\K$ is \emph{$d$-Matou{\v s}ek}, if the $\bbZ_2$-index of the space
$$\hat {\K}=\Set{(x,y)\in\|\K'\|^2\ :\ \Par{\bigcap\text{supp}(x)}\cap\Par{\bigcap\text{supp}(y)}\notin \K}$$
is less than $d$. Here $\text{supp}(x)$ denotes the {\sl support} of $x$ in $\K'$, which is the inclusion-minimal 
face of $\K'$ containing $x$. 

It is straightforward to verify that $\K$ is $d$-representable iff there exists a linear map $f:\K'\to\bbR^d$,
such that for every set $I\su V(\K)$ not in $\K$, we have $\bigcap_{i\in I}f(\alpha_i)=\emptyset$,
where $\alpha_i=\Set{J\in V(\K')\ :\ i\in J}$. This implies
the existence of a $\bbZ_2$-map from $\hat {\K}$ to $S^{d-1}$;
thus any $d$-representable complex is also $d$-Matou{\v s}ek.
White proved that nerves of good covers in $\mathbb R^d$ are $d$-Matou{\v s}ek
and also showed that being $1$-Matou{\v s}ek is equivalent to being $1$-representable.

\smallskip
Regarding part (iii), Debrunner~\cite{Debrunner:1970th} showed that for the statement of topological Helly it suffices to assume
that the (reduced) homology
of intersections of $k$ sets in the family $1 \le k \le d+1$  vanishes at and below dimension $d-k$, and even more general 
conditions were found by Montejano \cite{Montejano:2014ii}.

\subsection {Conditions for the fractional Helly property}
\label {s:hel-fra}

\begin {prob}
(i) Find geometric, topological, and combinatorial conditions that imply the fractional Helly property.

(i) Find geometric, topological, and combinatorial conditions that imply the strong fractional Helly property.
\end {prob}

We will mention here two conjectures regarding the fractional Helly property and two related theorems.
A class of simplicial complexes is {\sl hereditary} if it is closed under induced subcomplexes. Recall that 
for a simplicial complex $\K$, $f_i(\K)$ is the number of $i$-faces of $\K$ and $b(\K)$ is the sum of (reduced)
Betti numbers of $\K$.
In connection with the fractional Helly theorem, Kalai and Meshulam~\cite{kalai:2010}
formulated the following conjecture.

\begin {conj}[Kalai and Meshulam] 
\label {c:km4}
Let $C>0$ be a positive number. 
Let $\F$ be the hereditary family of simplicial complexes defined by the property that 
for every simplicial complex $\K \in {\F}$ with n vertices,
$$b(\K) \le C n^d.$$
Then for every $\alpha >0$ there is $\beta = \beta(d,C) >0$ such that  
$\K \in {\F}$ and $f_d(\K) \ge \alpha {{n \choose {d+1}}}$ imply $\dim (\K) \ge \beta n$.  
\end {conj}
The conclusion of the conjecture is referred to as the {\sl fractional Helly property of degree} $d$. 
Kalai and Meshulam further conjectured that the conclusion holds even if one replaces $b(\K)$ with $|\chi (\K)|$,
where $\chi (\K)$ is the Euler characteristic of $\K$.
When $d=0$, this conjecture is about graphs and it was proved (in its strong version)
in \cite{ChSSS:2020}; see also \cite{ScottSey:2020}.


\begin {conj}[Kalai and Meshulam]
\label {c:km5}
Let $\cal U$ be a family of sets in $\R ^d$. Suppose that for every intersection $L$ of $m$ members of $\cal U$, 
$b(L) \le \gamma m^{d+1}$. Then $\cal U$ satisfies the fractional Helly property of order $d$.
\end {conj}

In some special cases the fractional Helly property has been established.
For instance, Matou{\v s}ek~\cite{Matousek:2004cs} 
showed that families of sets with a bounded VC dimension in $\R^d$ satisfy the fractional
Helly property of order $d$. Another case is 
the so-called {\sl convex lattice sets}. These are sets of the form $\mathbb Z^d \cap C$
where $\mathbb Z^d$ is the lattice of integer points in $\R^d$ and 
$C$ is a convex sets in $\R^d$. A result of B\'ar\'any and Matou{\v s}ek~\cite{Barany:2003wg} asserts that 
families of convex lattice sets  in $\R^d$ satisfy the fractional Helly property of order $d+1$. 
In both of these theorems the fractional Helly number is considerably smaller than the Helly number.
For example, let $\F$ be the family of all convex lattice sets in $\R^d$. The Helly number of $\F$, ${\rm h}(\F)$, is equal to $2^d$, as shown by Doignon~\cite{Doignon:1973ht}, while the fractional Helly number is $d+1$; \cite{Barany:2003wg}.

\begin {prob}
  \label {p:rag}
Does the assertion of the Radon theorem imply the fractional Helly property? 
\end {prob}

An affirmative solution to one interpretation of this question was recently given by Holmsen and Lee~\cite {Holm:2019uv},
who showed that for abstract convexity spaces, the finite Radon number $r$ implies that the 
fractional Helly number is bounded by some function $m(r)$ of $r$.

\begin {prob}
  \label {p:rag2}
Estimate $m(r)$. 
\end {prob}

Convex sets are sets of solutions of systems of linear inequalities,
and we can consider systems of polynomial inequalities of higher degrees.

\begin {conj}
  The family ${\cal B}^d_k$
  of sets of solutions in $\mathbb R^d$ of polynomial inequalities of degree $\le k$
  has the fractional Helly property.
  \end {conj}

It is known \cite{Motzkin:1955} (and is an easy consequence of Helly's theorem itself)
that the class ${\cal A}^d_k$ of sets in $\mathbb R^d$ of common zeroes of systems of  
polynomial inequalities of degree $\le k$ has the Helly number ${{d+k} \choose {k}}$. And we can even ask
if this formula gives the precise fractional Helly number for the class ${\cal B}^d_k$. 

We conclude this section by mentioning an interesting recent abstract notion of convexity
described by Moran and Yehudayoff \cite {MorYe:2020},
which seems relevant to various
problems raised in this paper and, in particular, to Problem \ref{p:rag}.
In this notion of abstract convexity, which we call {\it MY-convexity},
we assume that every ``convex set'' is the
intersection of ``halfspaces.'' We assume further that the VC dimension of the class of halfspaces is at most $D$.
The class ${\cal B}_k^d$ is an example of an MY-convexity space where the hafspaces are the sets of
solutions of a single polynomial inequality of degree $k$. 

\begin {prob}
  Consider an MY-convexity space $X$ where the VC dimension of the class of halfspaces is at most $D$.
  (i) Does $X$ have the fractional Helly number $f(D)$ for some function $f$ of $D$?
  (ii) Does $X$ have the fractional Helly number $D$?
  \end {prob}


\subsection {The $(p,q)$-property}

The conclusion of Helly's theorem is that the family is {\it intersecting}; i.e.,
there is a point $\R^d$ that is included in all sets in the family.

\begin {prob}
What conditions guarantee that the family is $t$-pierceable, meaning that there are
$t$ points such that every set in the family contains at least one?

In the language of nerves, what conditions guarantee
that the set of vertices of the nerve can be expressed as the union of $t$ faces? 
\end {prob}

A family of sets has the $(p,q)$-{\sl property}
if for every $p$ members of the family some $q$ have a non-empty intersection. Note that here we assume $p\ge q > d$.
(For nerves this says that every set of $p$ vertices spans a face with $q$ vertices,
and this is closely related to Tur\'an's problem for hypergraphs.)

Hadwiger and Debrunner \cite {Hadwiger:1957we} introduced the $(p,q)$-property and
proved
\begin {theo} If a finite family of convex sets in $\R^d$ has the $(p,q)$-property and $(d-1)p < (q-1)d$, then it is $p-q+1$-pierceable.
\end {theo}
 
A family of sets has the $(p,q)_r$ {\sl property} if in its nerve every $p$ vertices span at least $r$ faces with $q$ vertices. 
This was introduced by Montejano and Sober\'on
\cite{Montejano:2011cz} and further studied by Keller and Smorodinsky~\cite{KellSmor:2018}. Montejano and Sober\'on
proved (among other results) 

\begin {theo} A family of convex sets in $\mathbb R^d$ with the $(d+2,d+1)_d$ property is 2-pierceable.  
\end {theo}

Hadwiger and Debrunner~\cite {Hadwiger:1957we} conjectured in 1957 and Alon and Kleitman~\cite{Alon:1992ta,Alon:1992gb} 
proved the following important theorem.

\begin {theo}[$(p,q)$-theorem] For all $p\ge q >d$ there exists $f(d,p,q)$ such that 
if a family of convex sets in $\R^d$ has the $(p,q)$-property, then it is $f(d,p,q)$-pierceable.
\end {theo}
 
The bound on $f(d,p,q)$ given in~\cite{Alon:1992ta} is enormous.
The first open case is $d=2$ and $p=4,q=3$. It is known that 
$f(2,4,3)$ is between $3$ and $9$, the lower bound is from~\cite{Kleitman:2001vo},
and the upper bound is a recent result of 
McGinnis~\cite{McGinn:2020} who brought down the upper
bound of $13$ of~\cite{Kleitman:2001vo} to $9$. Substantial improvements
for the general case were given by Keller, Smorodinsky, and Tardos~\cite{KellSmoTar:2018} and by Keller and Smorodinsky~\cite{KellSmor:2020}. 

\begin {prob} Improve further the bounds on $f(2,4,3)$ and, more generally, on $f(d,p,q)$.
\end {prob}

Alon, Kalai, Matou{\v s}ek, Meshulam~\cite{Alon:2002wz} proved the following result that implies that the 
Alon--Kleitman theorem extends
to good covers and Leray complexes (but with worse bounds).

\begin {theo}
\label {t:akmm}
For every $q > d+1 $ there exists $C(d,q)$ with the following property: 
let $\F$ be a hereditary class of simplicial complexes satisfying the fractional Helly property of degree $d.$ 
If a simplicial complex $\K \in \F$ has the property that every $q$ vertices span a $d$-dimensional face,   
then the vertices of $\K$ can be covered by $C(d,q)$ faces.
\end {theo}

See Eckhoff \cite {Eckhoff:2003ed} for a survey on $(p,q)$-theorems.

\subsection {A Ramsey type question}

\begin {conj}
  For integers $d \ge 1$ and $r > 1+\lceil d/2\rceil$, there is $\alpha=\alpha (d,r)>0$ such that the following holds:  
  let $\F$ be a family of $n$ convex sets in $\R^d$. Then  $\F$ contains $n^{\alpha (d,r)}$ sets
    such that either every $r$ has a point in common
    or no $r$ has a point in common.
\end {conj}

There is a large literature on this and related questions starting with a theorem of
Larman, Matou{\v s}ek, Pach, and T\"or\H ocsik~\cite{Larman:1994wt}
that proves the case $d=2$ and $r=2$. Subsequent works are~\cite{Alon:2005hj} and \cite{Fox:2011ca}.
When $r=d+1$ this conjecture holds with $\alpha=1/(d+1)$.
This was observed by Keller and Smorodinsky (private communication)
and follows from their improved $(p,q)$-theorems. 
The general phenomenon here (with several interesting manifestations) 
is that graphs and hypergraphs arising in geometry satisfy much stronger
forms of Ramsey's theorem than arbitrary graphs and hypergraphs.

\subsection {Colorful, fractional colorful, and matroidal Helly theorems}

The colorful Helly theorem of Lov\'asz (see~\cite{Barany:1982va}) asserts the following. Assume 
that $\C_1,\ldots,\C_{d+1}$ are finite families of convex sets in $\R^d$ 
with the property that every transversal $K_1,\ldots,K_{d+1}$ is intersecting, then $\bigcap \C_i \ne \emptyset$ for some $i \in [d+1]$. Here {\sl transversal} means that 
$K_i \in \C_i$ for every $i \in [d+1]$. The colorful version implies the original one when $\C_1=\ldots =C_{d+1}$. 

The analogous colorful version of the fractional Helly theorem says that if an $\alpha$ fraction of all transversals of the system $\C_1,\ldots,\C_{d+1}$ is intersecting, 
then one of the families, say $\C_i$, contains an intersecting subfamily of size $\beta|\C_i|$. Here $\alpha>0$, of course, and $\beta=\beta(d,\alpha)$ has to be positive. 
Such a theorem (with $\beta=\alpha/(d+1)$) was proved and used first in~\cite{Arocha:2009ft}. The dependence of $\beta$ was improved by Kim~\cite{Kim:2017}, 
who showed in particular that $\beta \to 1$ as $\alpha \to 1$. The optimal dependence of $\beta$ on $\alpha$ and $d$ is a recent result of Bulavka, Goodarzi, 
and Tancer~\cite{BuGooTan:2020}. They use Kalai's algebraic shifting technique~\cite{Kalai:1984bg} and raise the following interesting conjecture.
 
\begin{conj} Let $\K$ be a $d$-Leray simplicial complex whose vertex set $V$ is
  partitioned into sets $V_1,\ldots,V_{d+1}$, called colors, and $|V_i|=n_i$ for $i \in [d+1]$. Assume that
$\K$ contains at least $\alpha \prod_1^{d+1}n_i$ colorful $d$-faces for some $\alpha >0$. Then there is $i \in [d+1]$ such that the dimension of the restriction of $\K$ to $V_i$ 
is at least $(1-(1-\alpha)^{1/(d+1)})n_i-1$. 
\end{conj}

Kalai and Meshulam \cite {Kalai:2005cm} extended the assertion of the colorful Helly theorem
to the topological setting and also considered a matroidal version. A matroidal complex is the complex
consisting of the independent sets of a matroid. Equivalently, $M$ is a matroidal complex
if and only if every induced subcomplex is pure, i.e., if all
its maximal faces have the same cardinality.

\begin {theo}
\label {t:mat-hel}
  Let $X$ be a $d$-Leray complex on the vertex set $V$. Suppose that $M$ is a matroidal complex on the same
  vertex set $V$  with rank function $\rho$. If $M \subset X$, then there exists $\tau \in X$ with $\rho (V \backslash \tau) \le d$.
  
\end {theo}   

This theorem gives the colorful Helly property when $M$ is a transversal matroid
and it suggests a general way to extend
results about colorings. We will encounter
this idea in Section \ref {s:col-tve} below where we try to move from colorful versions of
Tverberg's theorem to matroidal versions.  
Theorem \ref{t:mat-hel} has interesting connections with advances in topological combinatorics
related to Hall's marriage theorem and ``rainbow'' matchings; see \cite{Aharoni:2006} and \cite{Aharoni:2009}.
  
\section {More around Helly's theorem}
\label {s:h2}
\subsection {Dimensions of intersections: Katchalski's theorems} 


Let $g(d,k)$ be the smallest integer with the following property: for every family of $n$ convex sets in $\mathbb R^d$, $n \ge g(d,k)$, such that
the dimension of intersection of every $g(d,k)$ sets in the family is at least $k$,  the dimension of intersection of
all members of the family is at least $k$. Helly's theorem asserts that $g(d,0)=d+1$.
In 1971 Katchalski~\cite{Katchalski:1971to} proved the following interesting result.

\begin {theo} $g(d,0)=d+1, g(d,k)=\max \{d+1, 2(d-k+1)\}$ if $1 \le k \le d$.
\end {theo}

Given a family ${\K}=\{K_1,K_2,\dots, K_n\}$ of convex sets in $\R^d$ and $J \subset [n]$, set $K(J)=\bigcap_{j \in J}K_j$ and write $d(J)=\dim K(J)$.
A further remarkable result of Katchalski~\cite{Katchalski:1978to} ``reconstructs'' the dimension of the intersection: 

\begin {theo} \label {t:kd2} Let ${\K}=\{K_1,K_2,\dots, K_n\}$ and ${\K'}=\{K'_1,K'_2,\dots, K'_n\}$
be two families of compact convex sets in $R^d$. If $d_{\K}(J)=d_{\K'}(J)$ for every $J$, $|J| \le d+1$, then
$d_{\K}(J)=d_{\K'}(J)$ for every $J$. 
\end {theo}

Katchalski actually proved a stronger statement,
namely, that the condition $d_{\K}(J)=d_{\K'}(J)$ for every $J$ with $(d+1)-\lfloor d/2 \rfloor \le |J| \le d+1$ 
suffices for the conclusion of Theorem \ref {t:kd2}. More generally he proved that for every $r \ge 1$ if 
$d_{\K}(J)=d_{\K'}(J)$ for every $J$ with $(d+r)-\lfloor d/(r+1) \rfloor \le |J| \le d+r$, then $d_{\K}(J)=d_{\K'}(J)$ for every $J$.

Define the $D$-nerve of a finite set of convex sets as its nerve $K$ where every face $S \in K$ is labeled by 
the dimension of $\bigcap_{i \in S}K_i$.
We can regard the $D$-nerve as a nested collection of simplicial complexes that correspond to intersections of dimension $\ge j$.  

\begin {prob}
  Explore combinatorial and topological properties of $D$-nerves of families of compact convex sets in $\mathbb R^d$.
\end {prob}

\subsection {Helly with volume}

Theorems about volumes of intersections are closely related to theorems about dimensions of
intersections. The natural question is, given a finite family $\F$ of convex sets in $\R^d$, what condition guarantees that 
the intersection $\bigcap \F$  not only is non-empty but also has volume at least one, say. The first result in this direction is in~\cite{Barany:1982ga}
of B\'ar\'any, Katchalski, and Pach.

\begin {theo}[Helly with volume]\label{t:volumes} Assume that $\F$ is a finite family of convex sets in $\R^d$, $|\F|\ge 2d$, such that the intersection 
of any $2d$ sets from $\F$ has volume at least one. Then $\vol (\bigcap \F) \ge d^{-2d^2}$.
\end{theo}

The example of the $2d$ halfspaces in $\R^d$ whose intersection is the unit cube shows
that the number $2d$ is the best possible in this theorem. 
In other words, $2d$ is the {\sl Helly number for volumes}. However, the bound  $d^{-2d^2}$  is not sharp and was improved 
first by Nasz\'odi~\cite{Naszodi:2015vi} to $(cd)^{-2d}$ and later by Brazitikos~\cite{Brazit:2017} to $(cd)^{-1.5d}.$
In both estimates, $c>0$ is a universal constant. The following conjecture is still open.

\begin {prob}\label{p:volum} Show that under the conditions of Theorem~\ref{t:volumes}
$\vol (\bigcap \F) \ge (cd)^{-d/2}$ where $c>0$ is a constant. 
\end {prob}

A similar result was established in~\cite{Barany:1982ga} for the diameter of the intersection. The Helly number is again $2d$. 
So if the intersection of any $2d$ sets from the family $\F$ has diameter at least one, then $\diam\;\F  \ge cd^{-d/2}$. 
This lower bound was improved in a series of recent papers: first by Brazitikos~\cite{Brazit:2017} to $cd^{-11/2}$, 
then by Ivanov and Nasz\'odi~\cite{NaIv:2021} to $(2d)^{-3},$ and most recently by Almendra-Hern\'andez, Ambrus, and Kendall~\cite{AlAmK:2021} 
to $(2d)^{-2}.$ This leads to the next problem.

\begin {prob}\label{p:diam}  Assume that $\F$ is a finite family of convex sets in $\R^d$, $|\F|\ge 2d$, such that the intersection 
of any $2d$ sets from $\F$ has diameter at least one. Then $\diam\; \bigcap \F \ge cd^{-1/2}$.
\end{prob}

Recently, several further quantitative Helly-type results have appeared; see for instance \cite{DFN:2020} and \cite{DillSob:2020}.

\subsection {Unions of convex sets: Around the Gr\"unbaum--Motzkin conjecture}

\label{s:h-u}

Nina Amenta~\cite{Amenta:1994gs} proved a Helly-type result on unions of disjoint convex sets.

\begin {theo}\label{t:amenta} Let $\F$ be a family of sets in $\R^d$ such that every member in
${\F}$ is the union of $k$ disjoint compact convex sets. Suppose further that every intersection of members of ${\F}$ is
also a union of $k$ disjoint convex sets. If every $k(d+1)$ sets in ${\F}$ has a point in common, then
$\bigcap {\F} \ne \emptyset$.
\end {theo}

In the language of Section~\ref{s:hnho}, Theorem~\ref{t:amenta} asserts that the Helly order of the family of
disjoint unions of $k$ compact convex sets in $\R^d$ is $(d+1)k$. This was conjectured by Gr\"unbaum and
Motzkin~\cite{Grunbaum:1961fd} who proved the case $k=2$; Larman~\cite{Larman:1968} proved their conjecture for $k=3$ and
Amenta in its full generality. It is easy to see that this family has no finite Helly number.  

Kalai and Meshulam \cite{Kalai:2008kc} proved that Amenta's theorem extends topologically.
They consider the following setting. Let $K$ and $L$ be simplicial complexes with a map from $V(K)$ to $V(L)$
such that the inverse image of every $i$ face in $L$ is the union of at most $k$ $i$-faces of $L$. If $K$
is $d$-Leray, then the Leray number of $L$ is at most $dk+k-1$. 

Eckhoff and Nischke \cite {Eckhoff:2009kv} showed that Amenta's theorem extends combinatorially. In the setting
of the previous paragraph they proved that if $K$ has no missing face of size $d+1$ or larger, then
$L$ has no missing face of size $k(d+1)$ or larger.

\subsection {More on families of unions of convex sets}

We may consider sets in $\R^d$ that can be represented as unions of $k$ convex sets but delete the disjointness
assumption. In this case Alon and Kalai ~\cite{Alon:1995fs} and Matou{\v s}ek~\cite {Matousek:1997di}
proved the following result.

\begin {theo} \label {t:akm} Assume that $\F$ is a finite family of sets in $\R^d$ such that every member in
$\F$ is the union of $k$  compact convex sets. Then $\F$ has a finite Helly order.
\end {theo} 

Let us mention a recent topological Helly-type theorem by Goaoc,  
Pat\'ak, Pat\'akov\'a, Tancer, and Wagner \cite{GPPTW:2017dp} that strengthens Theorem \ref {t:akm}.

\begin {theo}
\label {t:ter}
For every $\gamma >0$ there is $h(\gamma,d)$ with 
the following property: let $\cal U$ be a family of sets in $\mathbb R ^d$. 
Suppose that for every intersection $L$ of some members of $\cal U$ 
and every $i \le \lceil d/2 \rceil-1$, we have $b_i(L) \le \gamma $. 
Then, if every $h(\gamma,d)$ members of $\cal U$ have a point 
in common, then all sets in $\cal U$ have a point in common. 
\end {theo}

We note that Theorem \ref{t:ter} implies Theorem \ref{t:akm}. In fact, its proof relies on the method developed by Matou{\v s}ek
in \cite {Matousek:1997di}. His method, connecting topological obstructions for embeddability to Helly-type
theorems, is the basis of White's notion~\cite{White2020} of $d$-Matou{\v s}ek complexes.

\medskip
In connection with this we mention the following curious question.

\begin {conj}
  \label {c:curious}
  The Helly order of families of unions of two disjoint non-empty sets in $\mathbb R^d$ is $d+1$.
  \end {conj}


This is known to be false if ``two'' is replaced by a large integer even when $d=2$.

We say that two compact sets intersect nicely if the long Meyer--Vietoris exact sequence splits into short exact sequences
dimensionwise. 
\begin {prob}
  \label{nervg}
  Let ${\cal K} =\{K_1,K_2,\dots, K_n\}$ be a finite family of compact sets such that
    for every set of indices $I\subset [n]$,    ${\cal K}(I)$ is topologically equivalent to a fixed topological space $Z$, and
    for every two sets of indices $I,J \subset [n]$,    ${\cal K}(I)$
    and ${\cal K}(J)$ intersect nicely. Then $\bigcup {\cal K}$ is topologically equivalent to  a fiber bundle over ${\cal N}({\cal K})$
    with fibers topologically equivalent to $Z$. 
\end {prob}

A positive answer to Problem \ref {nervg} would imply Conjecture \ref {c:curious} because a pair of disjoint unions of non-empty convex sets
whose intersection is also a disjoint union of non-empty convex sets always intersect nicely.


\subsection {A conjecture by Gao, Landberg, and Schulman}

Here is an interesting Helly-type conjecture by Gao, Langberg, and Schulman~\cite{GaoLS:2008ej}.
For a convex set $K$ 
in $\R ^d$ an $\epsilon$ enlargement of $K$ is $K+\epsilon (K-K)$ (where $K-K=\{x-y: x, y \in K\}$).

\begin {conjecture}
\label {c:ter}
For every $d$, $k$, and $\epsilon$ there is some $h=h(d,k, \epsilon)$ with the following property.
Let $\F$ be a family of unions of $k$ convex sets. Let $\F^\epsilon$ be the family obtained by 
enlarging all the involved convex sets by $\epsilon$. If every $h$ members of $\F$ have a point in 
common, then all members of ${\cal F}^\epsilon$ have a point in common.
\end {conjecture}

Of course, for $k=1$ we can take $\epsilon =0$ by Helly's theorem.

\subsection {Boxes and products}

\begin {prob}
  Let $d_1,d_2,\dots, d_r$ be a partition of $d$. Study Helly-type theorems for families of
  Cartesian products $K_1 \times K_2 \times \ldots \times K_r$ of convex sets where
  $\dim K_i=d_i$.
\end {prob} 

The case of standard boxes, namely when $d_1=d_2=\cdots=d_{d}=1$ is of special interest.
Standard boxes have Helly number 2, and therefore
their nerves are determined by their graphs. Eckhoff proved an upper bound theorem for standard boxes \cite{Eckhoff:1988vk},
and studied the extremal families \cite {Eckhoff:1991ue}.
It is easiest to describe the
families where the upper bound is attained. If $f_{d+r}=0$ (that is, the largest non-empty intersection is for $d+r$ sets), then 
the family consists of $r$ copies of $\mathbb R^d$ and roughly the same number of parallel copies of
each of the $d$ coordinate's hyperplanes.  

Let $\K$ be the nerve of a family of standard boxes in $\mathbb R^d$. Then $\K$ is a $d$-Leray complex and has 
the further property that if $S$ is a set of vertices such that every pair of vertices in $S$ form an edge, then $S$ is a face of $\K$. 
This property of the nerve corresponds to Helly number 2 for the original family and we refer to it as Helly number 2. 

  \begin {prob}
    Extend Eckhoff's upper bound theorem to the class of  $d$-Leray complexes with no missing faces of size
    greater than 3 (namely, those corresponding to Helly number 2).
  \end {prob}

\subsection {Mutual position of convex sets}

The study of nerves of convex sets is the study of intersection patterns of families of convex sets.
When we start with a family of convex sets in $\R^d$ we can go further and consider intersection patterns
of the convex hulls of all subfamilies. (We can go even further by alternating between taking
convex hulls and intersections and by considering statements regarding $k$-flat transversals rather than plain intersections.) 

\begin{figure}[h!]
\centering
\includegraphics[scale=0.7]{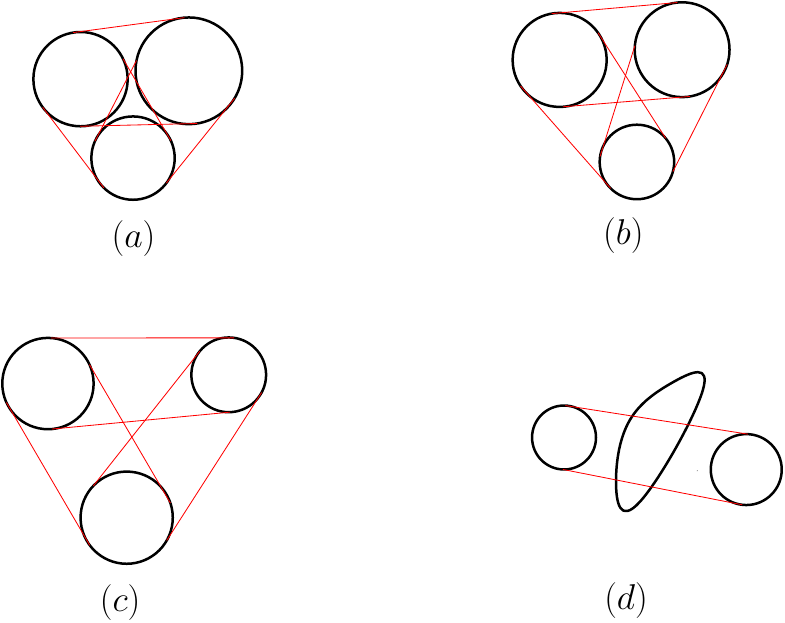}
\caption{Mutual positions of three convex sets.}
\label{fig:gil}
\end{figure}
Figure~\ref{fig:gil} shows various possible positions of three convex sets in the plane:\\
\indent (a) the convex hull of every two sets intersects the third set,\\
\indent (b) the convex hull of any two sets is disjoint from the third set, but 
all pairwise convex hulls have a point in common,\\
\indent (c) the three convex hulls of pairs of sets have no point in common,\\
\indent (d) the convex hull of two sets intersect the third set.

\smallskip
Statements in this wider language can be regarded as the study
of mutual positions of convex sets and they are, of course, of interest even for configurations of points,
which we discuss in the next sections.

\begin {prob}\label{mutual}
  Are there interesting things to say about the mutual position of convex sets?
  \end {prob}

\subsection {Order types for points and sets}

To conclude this section and prepare for the next, we briefly mention the notion of order types (aka oriented matroids).
These objects arise from configurations of points (or of hyperplanes) in real vector spaces,
and can also be associated with directed graphs.
Consider a sequence $Y=(y_1,y_2,\dots,y_n)$  of $n$ points in $\mathbb R^d$ that affinely span $\mathbb R^d$.
The order type 
described by $Y$ can be seen as the set of all minimal Radon partitions. There is a more general axiomatic
definition of order types that roughly requires that the restriction to every $d+3$ points be 
an order type of $d+3$ points
in a real space. For general order types there is a topological representation that replaces the linear
description of order types that correspond to point configurations.
Another equivalent way to describe the order type is as follows:
for every set $J$ of subscripts $i_1,\ldots,i_{d+1}$ with $1\le i_1 < \ldots < i_{d+1} \le n$, let ${\rm sg}(J,Y)$
be the sign of the determinant of the
$(d+1)\times (d+1)$ matrix
\begin{equation}
\begin{pmatrix}
y_{i_1} & y_{i_2}& \cdots  & y_{i_{d+1}}\\
1          & 1        & \cdots  & 1
\end{pmatrix}
.
\end{equation}
Two sequences  $Y=(y_1,y_2,\dots,y_n)$ and $Z=(z_1,z_2,\dots,z_n)$  of $n$ points in $\mathbb R^d$ are {\sl equivalent} 
(or have the {\sl same order type}) if ${\rm sg}(J,Y)={\rm sg}(J,Z)$ for all $J\subset [n]$ of size $d+1$.  


For more on oriented matroids see \cite{Bjor:1993}.
Returning to families of convex sets we note that one way to record the mutual position of $n$ convex
sets $K_1,K_2,\dots,K_n$ in $\mathbb R^d$ is by listing all order types of sequences
$y_1 \in K_1, y_2 \in K_2, \dots, y_n \in K_n$. 

Goodman and Pollack's notion of allowable sequences for 
configurations ~\cite{Goodman:1993ji} is a very useful way to study order types of planar configurations. 
The more general notion of interval sequences by Dhandapani, Goodman, Holmsen, and Pollack 
gives a way to record mutual positions
of $n$ convex planar sets~\cite{Dhand:2005}.

\section {Around Tverberg's theorem}
\label{s:t1}

\subsection {Sierksma's conjecture}

\begin {conj} [Sierksma's]\label{sierksma} 
The number of Tverberg $r$-partitions of a set 
of $(r-1)(d+1)+1$ points in $ \mathbb R ^d$ is at least $ ((r-1)!)^d$.
\end {conj}

This question was raised by Sierksma~\cite{Sierksma:1979} and not much progress has been achieved since. The best lower bound is about the square root of the 
conjectured one. This is a result of \cite{Vucic:1993be} and \cite{Hell:2007tp}. The conjecture, if true, is sharp, as shown by the example in Figure 3 for $d=2,r=4$: 
the vertices of the 3 triangles plus the point in the center is a set with 10 points and $3!^2$ Tverberg partitions. 

\begin{figure}[h!]
\centering
\includegraphics[scale=0.6]{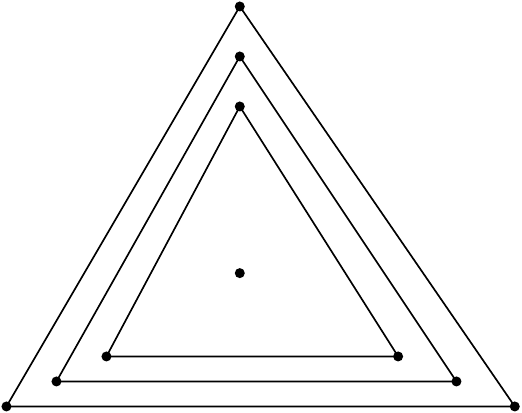}
\caption{10 points with $3!^2$ Tverberg partitions.}
\label{fig:Sierks}
\end{figure} 

In $\R^d$ take analogously $r-1$ 
$d$-dimensional simplices with their center at the origin; their vertices together with the origin form a set of $(r-1)(d+1)+1$ points with $(r-1)!^d$ Tverberg partitions. 
There are further cases where equality holds, such as the one connected to the following problem raised by Perles. We need a definition: a Tverberg partition $S_1,\ldots,S_r$ of an $m$-element set $X\subset \R^d$ is {\sl of type} $(a_1,a_2,\ldots, a_r)$ if the multisets $\{a_1,a_2,\ldots, a_r\}$ and  $\{|S_1|,|S_2|,\ldots, |S_r|\}$ coincide.

\begin {problem} 
Suppose that $a_1,a_2, \ldots, a_r$ is a partition of $m=(r-1)(d+1)+1$ such 
that $1\le a_i\le d+1$ for every $i$. Is there a configuration of $m$ points in $\mathbb R ^d$ for which all  
of Tverberg partitions are of type $(a_1,a_2,\ldots, a_r)$?  
\end {problem}

This problem was raised by Perles many years ago and a positive answer was given by White~\cite{Whi17}. 
White's examples provide a rich family of examples for cases of equality in Sierksma's conjecture. An even more general family of constructions 
for the equality cases, based on staircase convexity, is in the paper of Bukh, Loh, and Nivasch ~\cite{Bukh:2017vs}.
A similar construction was given by  P\'or~\cite{pornew} in connection with the so-called universal Tverberg partitions.

\begin {prob}
  Explore further examples of equality cases in Sierksma's conjecture.
  \end {prob}

\subsection  {Topological Tverberg}
\label {s:toptve}
 
\begin {conj} [topological Tverberg conjecture]
Let $ f$ be a continuous function from the $ m$-dimensional 
simplex $ \sigma^m$ to $ \mathbb R ^d$. If $ m \ge (d+1)(r-1)$
then there are $ r$ pairwise disjoint 
faces of $ \sigma^m$ whose images have a point in common.
\end {conj}

If $ f$ is a linear function, this conjecture reduces to Tverberg's theorem.
The case $ r=2$ was proved
by Bajm\'oczy and B\'ar\'any~\cite{Bajmoczy:1979bj} using the Borsuk--Ulam theorem. 
Moreover, for $r=2$, one can replace the simplex by any other 
polytope of the same dimension. 
The case where  $r$ is a prime number was proved in an important 
paper of B\'ar\'any, Shlosman, and  Sz\H ucs~\cite{Barany:1981vh}. 
The prime power case was settled by \"Ozaydin, in an unpublished (yet available) paper~\cite{Oza87}. 
For the prime power case,  the proofs are quite difficult and are based on computations of certain characteristic classes.

In 2015 the topological Tverberg conjecture was disproved in a short note by Frick~\cite{Frick:2015}. This involves 
some early result on vanishing of topological obstructions by \"Ozaydin, 
a theory developed by Mabillard and Wagner~\cite{MW14SoCG} extending Whitney's trick 
to $k$-fold intersections, and a fruitful reduction by Gromov~\cite{Gromov:2010eb}, rediscovered and extended by   
Blagojevi\'c, Frick, and Ziegler~\cite{Blagojevic:2017js}.  

\begin {conj} 
Let $ f$ be a linear function from an $ m$-dimensional 
polytope $P$ to $ \mathbb R ^d$. If $ m \ge (d+1)(r-1)$, then there 
are $ r$ pairwise disjoint faces of $P$ whose images have a point in common.
\end {conj}

\begin {prob} 
Does the conclusion of the topological Tverberg conjecture hold if the images of the faces under $f$ form a ``good cover''
(that is if all those images and all their non-empty intersections are contractible)?
\end {prob}

\subsection {Colorful Tverberg}
\label {s:col-tve}

Let $ C_1,\ldots,C_{d+1}$ be disjoint subsets of $ \mathbb R ^d$, called colors, each of 
cardinality at least $ t$. A $ (d+1)$-subset $ S$ of $\RR^d$ is said 
to be multicolored (or rainbow) if $|S\cap C_i|=1$ for $ i=1,\ldots,d+1$. Let $ r$ be an integer, and 
let $ T(r,d)$ denote the smallest value $ t$ such that for every collection of colors $ C_1,\ldots,C_{d+1}$ of 
size at least $ t$ each there exist $ r$ disjoint multicolored sets $ S_1,\ldots,S_r$ such 
that $ \bigcap^r_{i=1}\conv S_i\not=\emptyset$. The question of finiteness of $T(r,d)$ was raised in~\cite{BarFL:1992}
and proved there for the case $d=2$.

The general case was solved by an important theorem of \v{Z}ivaljevi\'c and Vre\'cica~\cite{ZivV:1992}. It asserts that $ T(r,d)\leq 2r-1$ if $ r$ is a prime, which implies that 
$T(r,d)\leq 4r-1$ for all $ r$ and $d$. This theorem is one of the highlights of discrete geometry and topological combinatorics. 
The only known proofs of this theorem rely on topological arguments although the statement is about convex hulls, partitions, and 
linear algebra. The following question is a challenge for convex geometers.

\begin {problem} Find a non-topological proof of the finiteness of $T(r,d)$.
\end {problem}

B\'ar\'any and Larman~\cite{Barany:1992tx} showed that $T(r,2)=r$ and asked the following.

\begin {conj}[colorful Tverberg conjecture] $ T(r,d)= r.$
\end {conj}

The case where $r+1$ is a prime was proved by Blagojevi\'c, Matschke, and Ziegler~\cite{blago15}. It is a neat result of Lov\'asz 
that appeared in~\cite{Barany:1992tx} that $T(2,d)=2$ for all $d$. Sober\'on gives an equally neat (and very different) proof of the same result in~\cite{Sober:2015}.

The colorful Tverberg theorem is related to a
well-known problem in discrete geometry, that of halving lines and hyperplanes. Given $2n$ points in general position in $\mathbb R^d$,
a halving hyperplane is a hyperplane with $n$ points on each side.

\begin {prob}
  What is the maximum number $H(2n,d)$ of partitions of a set of $2n$ points in $\mathbb R^d$ into equal parts via halving hyperplanes?
  Equivalently, what is the minimum number of non-Radon partitions with parts of equal size?
  \end {prob}
A well-known conjecture that is open even for $d=2$ is that for a fixed $d$, $H(n,d)=n^{d-1+o(1)}$.
With the help of the colorful Tverberg theorem it was shown that $H(n,d)=n^{d-\epsilon_d}$, where $\epsilon_d$ is a positive
constant depending on $d$. For $d=2$ it is known that $$n e^ {C \sqrt n} H(2,n) \le O(n^{4/3}).$$

A matroid version of Tverberg's theorem is the topic of~\cite{matroid-tverberg16}, which states 
the following.
Assume that $M$ is a matroid of rank $d+1$.
Let $b(M)$ denote the maximal number of disjoint bases in $M$.
If $f$ is a continuous map from the matroidal complex of $M$ to $\R^d$, then there exist $t\ge \frac 14 \sqrt {b(M)}$ 
independent sets $\sigma_1, \ldots, \sigma_t \in M$ such that $\bigcap_1^t f(\sigma_i)\ne \emptyset$.
It is not clear how good this lower estimate on $t$ is. 

\begin {conj}
  In the above theorem, $\sqrt {b(M)}$ could be replaced by $c b(M)$ for some absolute positive
  constant $c$.
  \end {conj}

\section {The cascade conjecture and more} 
\label {s:t2}

When we have $r<d+2$ points in $\mathbb R^d$ they have a Radon partition iff they are affinely dependent.
Are there conditions that guarantee that the existence of Tverberg partitions below the Tverberg number?
In this section we will discuss the dimension of Tverberg points and the quest for conditions guaranteeing the 
existence of Tverberg partitions for configurations of points below the Tverberg number. 

\subsection {The cascade conjecture} 

For a set $A$, denote by $T_r(A)$ the set of points in $\mathbb R ^d$ that belong to 
the convex hull of $r$ pairwise disjoint subsets of $ A$. We call these points Tverberg points of order $r$. 

Let $\bar t_r(A)= 1+ \dim T_r(A).$ (Note that $\dim\; \emptyset = -1$.) Radon's theorem can be stated as follows: if $\bar t_1(A) < |A|$ then $T_2(A) \ne \emptyset$. A similar statement
which is still open is: if $\bar t_1(A)+\bar t_2(A) < |A|$ then $T_3(A) \ne \emptyset$. We can go one step further: 
if $\bar t_1(A)+\bar t_2(A) +\bar t_3(A)< |A|$ then $T_4(A) \ne \emptyset$.
These statements are special cases of  

\begin {conjecture}[cascade conjecture]
\label{cascade}
For every $A \subset \mathbb R ^d$, $$\sum_{r \ge 1} \bar t_r(A) \ge |A|.$$
\end {conjecture}

This is a question of Kalai from 1974~\cite{kalai:1994}; see also \cite{Kal00flavor}.
The conjecture was proved for $d \le 2$ by Akiva Kadari (unpublished MSc thesis in Hebrew).
While this conjecture is wide open we can ask for topological extensions of various kinds and for 
more general topological conditions for configurations of cardinality below the Tverberg 
number $(r-1)(d+1)+1$,
that imply the existence of a Tverberg partition into $r$ parts; see Problem \ref{p:etp}.  

\subsection {Reay's dimension conjecture} 
\label {s:reay1}

The following is a 1979 question from Reay~\cite{Reay79} where general position means {\sl weak} general position;
that is, no $d+1$ points lie in a hyperplane.

\begin {conj} [Reay's conjecture] If $A$ is a set of $(d+1)(r-1)+1+k$ points in general position in $\mathbb R ^d$, then 
$${\rm dim}\; T_r(A) \ge k.$$
\end {conj}

In particular, Reay's conjecture asserts that a set of  $(d+1)r$ points in general position in $\mathbb R ^d$ can be partitioned 
into $r$ sets of size $d+1$ such that the simplices described by these sets have an interior point in common. This is easy when the points are 
in very general position, for instance, when they are algebraically independent. 
The main difficulty is how to use the weak general position condition. A recent result of Frick and Sober\'on~\cite{FriSob:2020}  
(see Section~\ref{s:c1}) is perhaps relevant here. 
While the conclusion of the cascade conjecture seems stronger than that of Reay's dimension conjecture, it is not known how to derive it from the cascade conjecture.

\subsection {Special cases of the cascade conjecture and expressing a directed graph as union of two trees}

A special case of the cascade conjecture  asserts that given $2d+2$ points in $\mathbb R ^d$, you can 
either partition them into two simplices whose interiors intersect, or 
you can find a Tverberg partition into three parts. A reformulation based on positive hulls is:
given $2d$ non-zero vectors in $\mathbb R ^d$ such that the origin is a vertex of the cone spanned by them, 
it is the case that either:

\begin {itemize}
\item

We can divide the points into two sets $A$ and $B$ so that the cones
 spanned by them have a $d$-dimensional intersection, or

\item
We can divide them into three sets $A$, $B$, and $C$ so that the cones 
spanned by them have a non-trivial intersection.

\end {itemize}

Another interesting reformulation is obtained when 
we dualize using the Gale transform, and this has led to the problem we consider next:
a very special class of configurations arising from graphs. 
Start from a directed graph $G$ with $n$ vertices
and $2n-2$ edges and associate with each directed edge $(i,j)$ 
the vector $e_j-e_i$. This leads to the following problem.

\vspace{0.1in}

\begin {problem}
Let $G$ be a directed graph with $n$ vertices and $2n-2$ edges. When can we divide the set of edges 
into two trees $T_1$ and $T_2$ (we disregard the orientation of edges) 
so that when we reverse the directions of all edges in $T_2$ we get a strongly connected digraph?
\end {problem}

One of us (Kalai) conjectured that if $G$ can be written as the union of two trees, 
the only additional obstruction is that there is a cut consisting only of two edges in reversed directions.
Chudnovsky and Seymour found an additional necessary condition: 
there is no induced cycle $v_1,v_k,\ldots,v_{2k},v_1$ in 
$G$, such that each vertex $c_i$ is cubic, the edges of the cycle alternate in direction, 
and none of the vertices $v_1,\ldots,v_{2k}$ are sources or sinks of $G$.

\subsection {Tverberg partitions of order three for configurations below the Tverberg number}
\label{s:Radon}

\begin{prob}\label{p:etp}
  When $n < 2d+3$, find conditions for the set of Radon points and the set of Radon partitions 
  of a set $X$ of $n$ points 
in $\mathbb R^d$, that guarantee the existence of a Tverberg partition into three parts.
 \end {prob}

The cascade conjecture asserts that if $n=d+2+k$ and the dimension of Radon
points is smaller than $k$, then there exists a Tverberg partition into three
parts. While this is wide open, it would be interesting to propose a more general topological
condition that suffices for the existence of a Tverberg partition into $r$ parts.

\begin {conj}
\label {vague}
If the map from the Radon partitions of $X$ to the Radon points of $X$ is topologically degenerate (in some sense), then a Tverberg partition into three parts exists. 
\end {conj}

In Problem~\ref{p:etp} and Conjecture~\ref{vague} we can relax the conclusion and can do so in various ways.   
For that we need a few definitions: the $k$-{\sl core} of a finite set $X$ in $\RR^d$ is 
the intersection of the convex hull of all sets $A\subset X$ with $|X\setminus A|\le k$,
that is, 
\[
{\rm core}_k\; X=\bigcap\{\conv A: A\subset X,\; |X\setminus A|\le k\}.
\] 
The case $k=0$ is the usual convex hull. 
The {\sl $k$-Radon core} of a finite set $X$ in $\RR^d$ is 
the intersection of Radon points of all sets $A\subset X$ with $|X\setminus A|\le k$; this is the
set of points in $\R^d$ that remain Radon points of $X$ even after we delete $k$ points from $X$ 
in all possible ways. (Clearly, the Tverberg points of order three are in the first Radon core,
and the points in the first Radon core are in the 2-core.)

\begin{prob}\label{p:etp2}
When $n < 2d+3$, find conditions for the set of Radon points and the set of Radon partitions 
  of a set $X$ of $n$ points in $\mathbb R^d$ that guarantee

(1) the second core of $X$, ${\rm core}_2\; X$, is non-empty,

(2) the first Radon core of $X$ is non-empty, 

(3) $X$ admits a Reay (3,2)-partition, that is, a partition into three parts such that the convex hulls
are pairwise intersecting; see Section~\ref{s:Reay}.
\end {prob}

\subsection {Radon partitions and Radon points for configurations based on cubic graphs}\label{s:cubic}

Let $G$ be a cubic graph with $n$ vertices $\{v_1,v_2,\dots, v_n\}$.
Associate with every edge $\{v_i,v_j\}$ in $G$ its characteristic vector in $\R^d$, giving a
configuration ${\rm Conf} (G)$ of $3n/2$ points in $(n-1)$-dimensional space. In~\cite{Onn:2001}
and also in personal communication (2011), Onn observed that the
existence of a Tverberg 3-partition (or even of a Reay (3,2)-partition; see Section~\ref{s:Reay})
is equivalent to a 3-edge coloring for $G$, and concluded that
deciding if a configuration of  $3(d+1)/2$ points in $\R^d$ ($d$ an odd integer)
admits a Tverberg partition into three part is NP-complete.

The following problem is motivated by the four color theorem.
\begin {prob}
  \label {p:etpc}
  (i) Study Radon partitions and Radon points for configurations based on cubic graphs.
  
   (ii) Find conditions for the Radon points and Radon partitions of\; ${\rm Conf}(G)$ that guarantee a 3-edge coloring for $G$.   
\end {prob}


It would be interesting to find conditions for Problems~\ref{p:etp} and Conjecture~\ref{vague}
that would imply the 3-edge colorability of bipartite cubic graphs
and, much more ambitiously, conditions that would imply the four-color theorem, namely, the 3-edge colorability of planar cubic graphs.

\section{More around Tverberg's theorem}
\label{s:t3}
\subsection  {Eckhoff's partition conjecture}
\label{s:tve-par}

Let $ X$  be a set endowed with an abstract closure operation $ X \to {\rm cl}(X)$. The only requirements of the closure operation are:

(1) $ {\rm cl}({\rm cl}(X))={\rm cl}(X)$ and

(2) $ A \subset B$ implies ${\rm cl}(A) \subset {\rm cl}(B)$.

Define $ t_r(X)$ to be the largest size of a (multi)set in $ X$ that cannot be partitioned into $r$ parts whose closures have a point in common. 
The following conjecture is due to Eckhoff~\cite{Eckhoff:2000jw}.

\begin {conj} [Eckhoff's partition conjecture] 
For every closure operation, $$ t_r \le t_2 \cdot (r-1).$$
\end {conj}

If $ X$ is the set of subsets of $ \mathbb R ^d$ and $ cl(A)$ is the convex hull operation, then Radon's 
theorem asserts that $ t_2(X)=d+1$ and Eckhoff's partition conjecture implies Tverberg's theorem. 
In 2010 Eckhoff's partition conjecture was refuted by Boris Bukh~\cite{bukh2010radon}. 
Bukh's beautiful paper contains several important ideas and further results. We will mention one ingredient. 
Recall the nerve construction for moving from a family $\F$ of $n$ convex sets 
to the simplicial complex that records empty and non-empty intersections for all subfamilies $\G$ of $\F$. 
Bukh studied simplicial complexes 
whose vertex sets correspond to the power set of a set of size $n$: 
starting with $n$ points in $\mathbb R ^d$ or 
some abstract convexity space, consider the nerve of convex hulls of all $2^n$ subsets of these points! 

In Bukh's counterexample, $ t_r = t_2 \cdot (r-1)+1$, which is just one larger than the conjectured bound. Perhaps $ t_r \le t_2 \cdot (r-1)+c$ for some universal constant $c\ge 1$.
There is a recent and positive development about Eckhoff's conjecture. P\'alv\"olgyi~\cite{Palv2020} has proved that $t_r$ grows linearly in $r$, 
that is, $t_r \le c r$ where the constant $c$ depends only on $r_2$.

\begin {prob}  
Find classes of closure operations for which  $$ t_r \le t_2 \cdot (r-1).$$
\end {prob}

We can ask if the inequality $ t_r \le t_2 \cdot (r-1)$ holds for  Moran and Yehudayoff's
 convexity spaces considered in Section \ref {s:hel-fra}.

Bukh's paper includes an interesting notion that extends the notion of nerves. 
Given a configuration of points in the Euclidean space or
in an abstract convexity space, we consider the nerve of convex hulls of all non-empty subsets of the points.
This is a simplicial complex that we refer to as the $B${\it-nerve} of the configuration, with the additional structure
that vertices are labeled by subsets, and with some additional combinatorial properties.

\begin {prob}  
  Study properties of $B$-nerves of point configurations in $\mathbb R^d$.
\end {prob}

\subsection {A conjecture by Reay}\label{s:Reay}

For a set $X \subset \mathbb R ^d$ a Reay $(r,j)$-{\sl partition} 
is a partition of $X$ into subsets $S_1,S_2,\ldots, S_r$ such that $\bigcap_{i=1}^j \conv S_{k_i} \ne \emptyset$,
for every $1 \le k_1 < \cdots < k_j \le r$. In other words, the convex hulls of any $j$ sets of the partition intersect. 
Define $R(d,r,j)$ as the smallest integer $m$ such that every $m$-element set $X \subset \mathbb R^d$ has a 
Reay $(r,j)$-partition. Reay~\cite{Reay79} conjectured that you cannot improve the value given by Tverberg's theorem, namely, that 

\begin {conj} [Reay's conjecture]
$R(d,r,j)= (r-1)(d+1)+1$. 
\end {conj}

Micha A. Perles believes that Reay's conjecture 
is false even for $j=2$ and $r=3$ for large dimensions, but with Moriah Sigron he proved~\cite{PeS16}
the strongest positive results in the direction of Reay's conjecture.

\subsection {Two old problems and universality}
\begin {problem}[McMullen and Larman]
  How many points $v(d)$ guarantee that for every set $X$ of $v(d)$ points in $\mathbb R^d$ there exists a 
  partition into two parts $X_1$ and $X_2$ such that for every $p \in X$,
  $$\conv (X_1 \backslash p) \cap \conv (X_2 \backslash p )\ne \emptyset.$$
\end {problem}

This is a strong form of Radon's theorem: the partition $X_1,X_2$ of $X=X_1\cup X_2$ remains a Radon partition even after we delete any point from $X$. Similar questions can be asked about Tverberg partitions. 
Larman~\cite{Larman:1072} proved that $v(d) \le 2d+3$ and this bound is sharp for $d=1,2,3,4$. The lower bound $v(d)\ge \lceil \frac {5d}3\rceil +3$ is a result of Ram\'{\i}rez Alfons\'{\i}n~\cite {Ram01}. This problem is the dual form of the original question by McMullen: what is the largest integer $n=f(d)$ such that every set of $n$ points in general position in 
$\mathbb R^d$ is projectively equivalent to the set of vertices of a convex polytope. 

\smallskip
A related problem is the following.

\begin {problem}
How many points $T(d;s,t)$ in $\mathbb R ^d$ guarantee that they can be divided into two parts such  
that every union of $s$ convex sets containing the first part has a non-empty intersection with every 
union of $t$ convex sets containing the second part.
\end {problem}

We explain next why $R(d; s,t)$ is finite. This is a fairly general Ramsey-type 
argument and it gives us an opportunity to mention a few recent important results.
The argument has two parts:

(1) Prove that $T(d;s,t)$ is finite (with good estimates) when the points are in cyclic 
position (to be defined shortly).

(2) Use the fact that for every $d$ and $n$ there is $f(d,n)$ such that 
among every $m$ points in general position in $\mathbb R ^d$, $m> f(d,n)$, one can find 
$n$ points in cyclic position.

\medskip

The finiteness of $T(d;r,s)$ follows (with horrible bounds) from these two ingredients by standard Ramsey-type results. 
It would be nice to understand the behavior of this function.

Statement (2) is a kind of {\sl universality theorem}. In a more precise form it says that for every $d$ and $n$ there is an integer $f(d,n)$
such that the following holds. Every sequence $x_1,\ldots,x_m$ in $\RR^d$ in general position with $m\ge f(d,n)$ contains a subsequence 
$y_1,\ldots,y_n$ such that all simplices of this subsequence are oriented the same way. The latter condition says, in more precise form, that 
for every set of subscripts $i_1,\ldots,i_{d+1}$ with $1\le i_1 < \ldots < i_{d+1} \le n$, the sign of the determinant of the $(d+1)\times (d+1)$matrix
\begin{equation}
\begin{pmatrix}
y_{i_1} & y_{i_2}& \cdots  & y_{i_{d+1}}\\
1          & 1        &  \cdots  & 1
\end{pmatrix}
\end{equation}
is the same (and different from $0$). Now a point set is {\sl cyclic} if its elements can be ordered so that the simplices along this ordering have the same orientation. 

Statement (2) says that the property of being cyclic is universal because every long enough sequence of points in general position contains a
cyclic subsequence of length $n$. Every finite sequence of points on the moment curve is cyclic. This shows that no other type 
of point sequence can be universal. Recently a fairly good understanding of $f(d,n)$ has been achieved in a series of papers.   
 
\begin {theo} $f(d,n) ={\rm twr}_d(\theta(n)).$
\end {theo}

Here, ${\rm twr}_d$ is the $d$-fold tower function. The lower bound is by Suk~\cite{Suk:2014tt} (improving 
earlier bounds by Conlon, Fox, Pach, Sudakov, and Suk~\cite{Conlon:2014fx}) and the upper bound comes from  
B\'ar\'any, Matou\v{s}ek, and P\'or~\cite{BarMP:2016}.

The following, somewhat vague, question emerges here naturally.

\begin{prob} Determine the universal type of $n$ lines in $\R^3$ and in $\R^d$. More generally, what is the universal type of $n$ 
$k$-dimensional affine flats in $\R^d$?
\end{prob}

Some preliminary results in this direction are the topic of a forthcoming paper by B\'ar\'any, Kalai, and P\'or~\cite{BarKP:2020}

We note that the order type of a sequence of points does not determine its Tverberg partitions.
\begin {prob}
  Develop a notion of order type based on Tverberg partitions into at most $r$ parts, $r \ge 3$.
\end {prob}

Here, Perles and Sigron's work on strong general position~\cite{PeS16}, and P\'or's
universality theorem~\cite{pornew} could be relevant.


\section{Carath\'eodory and weak $\epsilon$-nets} 
\label {s:c}

\subsection {Colorful Carath\'eodory and the Rota basis conjecture}\label{s:c1}

The following question was raised in Chow's Polymath~12 \cite{Polymath12} dedicated to Rota's basis conjecture.
Consider $d+1$ sets (or colors if you wish) $C_1,C_2,\cdots , C_{d+1}$ of points in $\mathbb R ^d$. 
Assume that each $|C_i|=d+1$ and that the interior of each $\conv C_i$ contains the origin.

\begin {prob}[D. H. J. Polymath]\label{polymath}
Can we find a partition of all points into $d+1$ rainbow parts such that the interior of the convex hulls of the parts 
have a point in common. (A rainbow set is a set containing one element from each $C_i$.) 
\end {prob}

To see the connection, first recall Rota's basis conjecture. 

\begin {conj}[Rota's basis conjecture]
If $B_1,B_2,\dots,B_n$ are disjoint bases in $\mathbb R^n$ (or even in an arbitrary matroid), then it is possible to find $n$ new disjoint bases $C_1,C_2,\dots,C_n$ such that each $C_i$ contains one element from every $B_j$. 
\end {conj}

Note that Rota's basis conjecture, can be stated (over $\mathbb R$) as follows:
Consider $d+1$ sets (or colors) $C_1,C_2,\cdots,C_{d+1}$ of points in $\mathbb R ^d$. 
Assume that each $|C_i|=d+1$ and that the interior of each $\conv C_i$ is nonempty. Then there exists a partition of all points into $d+1$ rainbow parts 
such that the interior of the convex hulls of each part is non-empty. 
 
Returning to Conjecture \ref {polymath}, we note here that according to the colorful Carath\'eodory theorem there is a rainbow set whose 
convex hull contains the origin. Without the words ``the interiors of'' Problem \ref {polymath} 
would be a special case of the colorful Tverberg conjecture (Section \ref {s:col-tve}).
A positive answer would be a  strong variant of Reay's conjecture (Section \ref{s:reay1}) 
on the dimension of Tverberg points,
and, as explained before, also a strong form of Rota's basis conjecture over the reals. 

A recent result of Frick and Sober\'on~\cite{FriSob:2020} is that
a set of $r(d+1)$ points in $\RR^d$ can always be partitioned into $r$ sets, each of size $d+1$, 
such that the convex hulls of the parts have a point in common. This theorem is related to the uncolored case of Problem~\ref{polymath}  
but without the word ``interior.''


\subsection {The complexity of the colorful Carath\'eodory theorem and of Tverberg partitions} 
\label {s:comp}
\begin {prob}
 Consider $d+1$ sets $C_1,C_2,\cdots , C_{d+1}$ of points in $\mathbb R ^d$. 
Assume that each $|C_i|=d+1$ and that each $\conv C_i$ contains the origin.
Is there a polynomial-time algorithm to find a rainbow simplex containing the origin? 
\end {prob} 

An interesting result in this direction is due to Meunier et al.~\cite{MeunM:2017}. They show that the problem lies in the intersection of 
complexity classes PPAD and PLS. The same applies to the analogous question about Tverberg partitions: is there a polynomial-time algorithm to find 
a Tverberg partition of an $(r-1)(d+1)+1$-element point set in $\R^d$? 
There are very few geometric problems in both classes PPAD and PLS that are not known 
to be solvable in polynomial time. The results in~\cite{MeunM:2017} are the first upper bound on the complexity of these problems.

\subsection {Carath\'eodory-type theorem for cores}

Recall the definition of the $k$-core of a finite set $X$ in $\R^d$ from Section~\ref{s:Radon}. The {\sl Carath\'eodory number} for the $k$-core is the smallest integer $f(d,k)$ with the property 
that $a \in {\rm core}_k X$ (where $X \subset \R^d$) implies the existence of $Y \subset X$ such that $a \in {\rm core}_k Y$ and $|Y| \le f(d,k)$. 
So $f(d,0)=d+1$ is just the Carath\'eodory theorem. B\'ar\'any and Perles~\cite{BarPer:90} established the finiteness of $f(d,k)$  together with 
some other properties of this function, for instance, that $f(d,1)=\max \{2(d+1),1+d+\lfloor d^2/4 \rfloor\}$, and that $f(2,k)=3(k+1)$. Several questions 
remain open; we mention only two of them.

\begin{prob} Determine $f(d,2)$ and $f(3,k)$.
\end{prob}

\subsection {The covering number theorem}

Assume that $X \subset \RR^d$ is finite and $|X|\ge d+1$. A {\sl simplex} of $X$ is just $\conv Y$ where $Y \subset X$ and $|Y|=d+1$.  
According to Carath\'eodory's theorem every point in $\conv X$ is contained in a simplex of $X$; that is, $\conv X$ is covered by the simplices of $X$. 
Which point is covered maximally, and how many times is it covered? A famous result of Boros and F\"uredi~\cite{Boros:1984ba}  says that 
in the planar case there is a point covered by $\frac 29 {n \choose 3}+O(n^2)$ simplices (that is, triangles) of $X$, where $n=|X|$. 
This is a positive fraction of all triangles of $X$ and the constant $\frac 29$ is the best possible. In higher dimensions Tverberg's theorem and the
colorful Carath\'eodory theorem imply (see \cite{Barany:1982va}) the following result.

\begin{theo}[covering number]\label{th:covnum}Assume $X$ is a set of $n\ge d+1$ points in $\RR^d$. 
Then there is a point covered by $\frac 1{(d+1)^d}{n \choose d+1}$ simplices of $X$.
\end{theo} 
This is again a positive fraction of all simplices of $X$. Define $b_d$ as the supremum of all $\beta>0$ such for that every set $X$ of $n\ge d+1$ points in $\RR^d$ 
there is a point covered by $\beta{n \choose d+1}$ simplices of $X$. So $b_d \ge (d+1)^{-d}$. 
In a remarkable paper, Gromov~\cite{Gromov:2010eb} showed that $b_d\ge \frac {2d}{(d+1)!(d+1)}$.
Gromov's theorem applies to continuous maps from the boundary of an $(n-1)$-dimensional simplex to $\mathbb R^d$. 
His estimate is an exponential improvement on the previous bounds. 
Both Gromov's theorem and Pach's theorem below play an important role in the emerging theory of high-dimensional expanders~\cite{FoxGr:2012}.

From the other direction Bukh, Matou\v{s}ek, and Nivasch~\cite{Bukh:2011vs} give an example, based on the stretched grid, 
that shows $b_d \le \frac {(d+1)!}{(d+1)^{d+1}}$. They conjecture that this is the right value of $b_d$.

\begin{conj} Show that $b_d = \frac {(d+1)!}{(d+1)^{d+1}}$. More modestly, prove that $b_d$ is exponential in $d$. 
\end{conj}

An interesting extension of the covering number theorem is the following result of Pach~\cite{Pach:1998vx}. 

\begin{theo}[Pach's theorem]\label{th:pach} Assume that $C_1,\ldots,C_{d+1}$ are sets (colors, if you like) in $\R^d$, each of size $n$. 
Then there is a point $p\in \R^d$ and there are subsets $D_i \subset C_i$ (for all $i\in [d+1]$), each of size at least $c(d)n$ 
such that the convex hull of every transversal of the system $D_1,\ldots,D_{d+1}$ contains $p$. Here $c(d)>0$ is a constant that depends only on $d$.
\end{theo}

This is a {\sl homogeneous} version of the covering number theorem. It was conjectured in~\cite{BarFL:1992}, 
where case $d=2$ was proved more generally even if the sets $C_1,C_2,C_3$ need not have the same size. This raises the following question.

\begin{prob} Does Pach's theorem remain true if the sets $C_1,\ldots,C_{d+1}$ have arbitrary sizes?
\end{prob}

We mention that Pach's theorem does not have a topological extension, as shown in~\cite{BarMesh:2018}, and in~\cite{BukHub:2020} in a stonger form.

\subsection {Weak $\epsilon$-nets}

An important application of the covering number theorem is about weak $\eps$-nets. Let $\eps>0$ be fixed. Given a finite set $X$ of $n\ge d+1$ points, let $\C$ 
be the (finite) family of sets $\conv Y$ for all $Y \subset X$ with $|Y|\ge \eps n$.
A set $F \subset \RR^d$ is called a {\sl weak} $\eps$-{\sl net} for $X$ if $F\cap C\ne \emptyset$ for every $C \in \C$. 

\begin{theo}[weak $\eps$-net theorem]\label{th:epsnet} Under the above conditions, there is a weak $\eps$-net $F$ for $X$ such that
\[ |F| \le \frac {c_d}{\eps^{d+1}},
\]
where $c_d>0$ is a constant. 
\end{theo}
The upper bound on the size of $F$ is from~\cite{Alon:1992ta} and\cite{Alon:1992ek} and has been improved to $O(\eps^{-d})$, disregarding 
some logarithmic terms. The trivial lower bound on the size of $F$ is $\frac 1{\eps}$. Bukh, Matou{\v s}ek, and Nivasch~\cite{Bukh:2011vs} 
give an example (based on the stretched grid or staircase convexity) where the size of the weak $\eps$-net is at least of order $\frac 1{\eps}(\log \frac 1{\eps})^{d-1}$. 
So the bounds on the size of a weak $\eps$-net are far from each other, and 
the general belief is that the true behavior should be slightly superlinear in $\frac {1}{\epsilon}$. 

\begin{prob} Find a better upper bound for the size of a weak $\eps$-net. 
\end{prob}

One remarkable improvement in this direction is a result of Rubin~\cite{Ru18} who showed that in the planar case 
there is always a weak $\eps$-net of size of order $\frac 1{\eps^{1.5+\delta}}$ for any $\delta >0$. 
A more recent result of Rubin~\cite{Ru20} applies in any dimension $d\ge 2$ and gives  a weak $\eps$-net of size 
of order $\frac 1{\eps^{d-1/2+\delta}}$ for any $\delta >0$.

\medskip
Weak $\eps$-nets can be defined not only for points but for $k$-dimensional affine flats in $\RR^d$. 
We only state the question for lines in $\RR^3$ and leave the rest of the cases to our imaginary reader. 
Let $\L$ be a set of $n$ lines and $\C$ be a finite family of convex sets in $\RR^3$. Assume 
that every $C \in \C$ intersects an $\eps$-fraction of the lines in $\L$, that is, 
\[
|\{L \in \L: C\cap L\ne \emptyset\}| \ge \eps n \mbox{ for every } C \in \C.
\] 

\begin{conj}[$\eps$-net of lines]\label{eps-net-lines} Under these conditions, there
  is a set of lines $\L^*$ whose size depends only on $\eps$ such that every $C\in \C$ intersects some line in $\L^*$.
\end{conj}

The set $\L^*$ can be thought of as a weak $\eps$-net of lines for $\C$.
We will encountered this question again soon, in connection with Problem~\ref{p:Rubin}.

\section {A glance at common transversals}
\label {s:ct}
\subsection {Transversals for intersecting families}

A $k$-transversal of a family of convex sets in $\R^d$ is a $k$-dimensional affine space that
intersects every set in the family. 
Transversal theory deals with conditions that guarantee the existence of $k$-transversals.
The case $k=0$ is connected to Helly-type theorems, and there are some general results for hyperplane transversals,
namely, $k=d-1$, and very few general results for $0<k<d-1$ and, in particular, for line transversals in $\mathbb R^3$.
The fascinating theory geometric transversals goes
beyond the scope of this paper; for surveys see Goodman, Pollack, Wenger~\cite {Goodman:1993xc}, Wenger~\cite{Wenger:2001}, 
and \cite{Holmsen:2013tu}. We will mention only a few problems where the conditions are in terms
of the intersection pattern of the sets in the family.

\begin{prob}\label{p:Rubin}
Assume that a family $\C$ of $n$ convex sets in $\R^3$ satisfies the the property that any two sets in $\C$ intersect. 
Show that there is a line intersecting $cn$ elements in $\C$, where $c>0$ is a universal constant.
\end{prob}

Partial results in this direction are given in~\cite{Barany:2021}.
Problem~\ref{p:Rubin} is the first, and so far most interesting, 
unsolved case of a series of problems of the same type. Namely, for what numbers $k,r,d$ is it true that, given a 
family $\C$ of convex sets in $\R^d$ where every $k$ tuple is intersecting, there is an $r$-flat intersecting a positive fraction of the sets in $\C$? 
Of course, the positive fraction should depend only on $k,r,$ and $d$.

\begin{figure}[h!]
\centering
\includegraphics[scale=0.6]{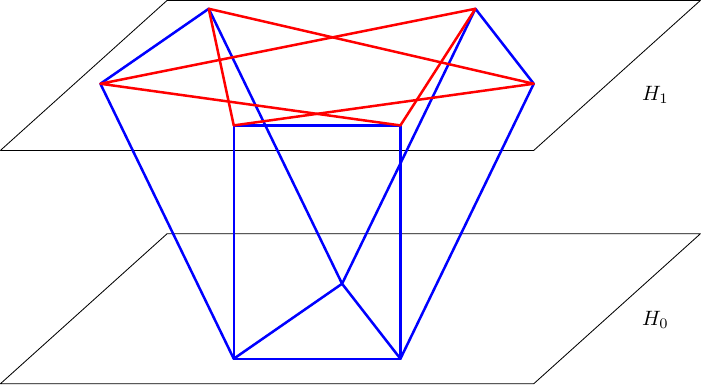}
\caption{Five sets in $\R^3$ that pairwise intersect and have no line transversal.}
\label{fig:rubin2}
\end{figure}

An interesting example satisfying the conditions is when $\C$ consists of $n$ lines in a two-dimensional plane in $\R^3$. 
Then, of course, every set in $\C$ is a line transversal for all sets in $\C$. This example shows that degenerate cases are going to make the problem difficult. 
Figure~\ref{fig:rubin2} is an example of five pairwise intersecting convex sets in $\R^3$ without a common line transversal. The five sets comprise three blue rectangles 
and two red triangles, all of whose vertices belong to two parallel planes $H_0$ and $H_1$.

\medskip

The question comes from a paper by Mart\'{\i}nez, Rold\'an, and Rubin~\cite{MarRR:2018} and is connected to the colorful Helly theorem. 
They also ask the slightly more general bipartite version of the question.

\vskip1pt 
\begin{prob}
Assume $\F$ and $\G$ are finite families of convex sets in $\R^3$ with the property that $A \cap B\ne \emptyset$ 
for any two sets $A \in \F$ and $B \in \G$. Show that there is a line intersecting $c|\F|$ elements of $\F$ 
or $c|\G|$ elements of $\G$ where $c>0$ is again a  universal constant.
\end{prob}
An example is two sets $\F$ and $\G$ of lines on a doubly ruled surface, which shows that degenerate cases may cause difficulties again. 
It is worth mentioning that both questions are invariant under non-degenerate affine transformation. 

We observe here that a positive answer to Conjecture~\ref{eps-net-lines} from the last section would imply that in Problem~\ref{p:Rubin}
there is a very finite set $\L$ of lines intersecting every element of $\C$, where by 
``very finite" we mean that the size of $\L$ is bounded by 1000, say, 
or by some other absolute constant.

\section {Conclusion}
\label {s:co}

This paper introduces the fascinating area of Helly-type theorems, and describes some of its main themes
and goals through a variety of open problems.  Often, results from convexity give a simple and strong
manifestation of theorems from topology: Helly's theorem manifests the nerve theorem from algebraic topology,
and Radon's theorem can be regarded as an early ``linear" appearance of the Borsuk--Ulam theorem.
One of our main themes is to further explore these connections to topology. Helly-type theorems also offer
complex and profound combinatorial connections and applications that represent the 
second main theme of this paper. We note that Helly-type theorems and the interplay between convex geometry,
combinatorics, and topology play an important role in the emerging theory of high-dimensional expanders.  

There are various related parts of this theory that we did not consider. We gave only a small
taste of the theory of common transversals; we did not discuss the closely related theorems of
Kirchberger and Krasnoselski{\v i}; and we did not consider the rich connections to metric
geometry. For example, when you consider families of translates of a fixed convex set
the theory takes interesting and surprising turns, and has applications and connections, e.g.,
to the theory of Banach spaces.

\bigskip
{\bf Acknowledgments.}  Research of IB was partially supported by Hungarian National Research grants (no. 131529, 131696, and 133819),
 and research of GK by the Israel Science Foundation (grant no. 1612/17).

\bibliographystyle{amsalpha}
\bibliography{references-2}

\newcommand{\etalchar}[1]{$^{#1}$}
\providecommand{\bysame}{\leavevmode\hbox to3em{\hrulefill}\thinspace}
\providecommand{\MR}{\relax\ifhmode\unskip\space\fi MR }
\providecommand{\MRhref}[2]{%
  \href{http://www.ams.org/mathscinet-getitem?mr=#1}{#2}
}
\providecommand{\href}[2]{#2}
\begin{thebibliography}{DLGMM19}

\bibitem[AB06]{Aharoni:2006}
R.~Aharoni and E.~Berger, \emph{{The intersection of a matroid and a simplicial
  complex}}, Trans. Amer. Math. Soc. \textbf{358} (2006), 253--267.

\bibitem[AB09]{Aharoni:2009}
\bysame, \emph{{Rainbow matchings in $r$-partite $r$-graphs}}, Electronic J.
  Combin. \textbf{16} (2009).

\bibitem[ABB{\etalchar{+}}09]{Arocha:2009ft}
J.~L. Arocha, I.~B{\'a}r{\'a}ny, J.~Bracho, R.~Fabila, and L.~Montejano,
  \emph{{Very colorful theorems}}, Discrete Comput. Geom. \textbf{42} (2009),
  142--154.

\bibitem[ABFK92]{Alon:1992ek}
N.~Alon, I.~B{\'a}r{\'a}ny, Z.~F{\"u}redi, and D.~J. Kleitman, \emph{{Point
  selections and weak $\varepsilon$-nets for convex hulls}}, Combin. Probab.
  Comput. \textbf{1} (1992), 189--200.

\bibitem[AHAK21]{AlAmK:2021}
V.~H. Almendra-Hern\'andez, G.~Ambrus, and M.~Kendall, \emph{Quantitative
  {H}elly-type theorems via sparse approximation}, 2021, arXiv:2108.05745,
  p.~6.

\bibitem[AK92a]{Alon:1992gb}
N.~Alon and D.~J. Kleitman, \emph{{Piercing convex sets}}, Bull. Amer. Math.
  Soc. \textbf{27} (1992), 252--256.

\bibitem[AK92b]{Alon:1992ta}
\bysame, \emph{{Piercing convex sets and the Hadwiger--Debrunner
  (p,q)-problem}}, Adv. Math. \textbf{96} (1992), 103--112.

\bibitem[AK95]{Alon:1995fs}
N.~Alon and G.~Kalai, \emph{{Bounding the piercing number}}, Discrete Comput.
  Geom. \textbf{13} (1995), 245--256.

\bibitem[AKMM02]{Alon:2002wz}
N.~Alon, G.~Kalai, J.~Matou{\v s}ek, and R.~Meshulam, \emph{{Transversal
  numbers for hypergraphs arising in geometry}}, Adv. in Appl. Math.
  \textbf{29} (2002), 79--101.

\bibitem[Ame94]{Amenta:1994gs}
N.~Amenta, \emph{{Helly-type theorems and generalized linear programming}},
  Discrete Comput. Geom. \textbf{12} (1994), 241--261.

\bibitem[APP{\etalchar{+}}05]{Alon:2005hj}
N.~Alon, J.~Pach, R.~Pinchasi, R.~Radoi{\v c}i{\'c}, and M.~Sharir,
  \emph{{Crossing patterns of semi-algebraic sets}}, J. Combin. Theory, Ser. A
  \textbf{111} (2005), 310--326.

\bibitem[B{\'{a}}r82]{Barany:1982va}
I.~B{\'{a}}r\'{a}ny, \emph{{A generalization of Carath{\'e}odory's theorem}},
  Discrete Math. \textbf{40} (1982), 141--152.

\bibitem[B{\'a}r21a]{Bar2021}
I.~B{\'a}r{\'a}ny, \emph{Combinatorial convexity}, AMS, 2021, in press.

\bibitem[B{\'a}r21b]{Barany:2021}
\bysame, \emph{{Pairwise intersecting convex bodies and cylinders in $\R^3$}},
  2021, arXiv:2104.02148, p.~9.

\bibitem[BB79]{Bajmoczy:1979bj}
E.~G. Bajm{\'o}czy and I.~B{\'a}r{\'a}ny, \emph{{On a common generalization of
  Borsuk's and Radon's theorem}}, Acta Math. Hung. \textbf{34} (1979),
  347--350.

\bibitem[BF84]{Boros:1984ba}
E.~Boros and Z.~F{\"u}redi, \emph{{The number of triangles covering the center
  of an n-set}}, Geom. Dedicata \textbf{17} (1984), 69--77.

\bibitem[BFL90]{BarFL:1992}
I.~B\'ar\'any, Z.~F\"uredi, and L.~Lov{\'a}sz, \emph{{On the number of halving
  planes}}, Combinatorica \textbf{10} (1990), 175--185.

\bibitem[BFZ19]{Blagojevic:2017js}
P.~V.~M. Blagojevi{\'c}, F.~Frick, and G.~M. Ziegler, \emph{{Barycenters of
  polytope skeleta and counterexamples to the topological Tverberg conjecture,
  via constraints}}, J. Europ. Math. Society \textbf{21} (2019), 2107--2116.

\bibitem[BGT20]{BuGooTan:2020}
A.~Bulavka, D.~Goodarzi, and M.~Tancer, \emph{{Optimal bounds for the colorful
  fractional Helly theorem}}, 2020, arXiv:2010.15765, p.~13.

\bibitem[BH20]{BukHub:2020}
B.~Bukh and A.~Hubard, \emph{On a topological version of {P}ach's overlap
  theorem}, Bull. London Math. Soc. \textbf{52} (2020), 275--282.

\bibitem[BKM17]{matroid-tverberg16}
I.~B{\'a}r{\'a}ny, G.~Kalai, and R.~Meshulam, \emph{A {Tverberg} type theorem
  for matroids}, Journey through Discrete Mathematics. A Tribute to
  Ji\v{r}\'{\i} Matou\v{s}ek (M.~Loebl, J.~Ne\v{s}et\v{r}il, and R.~Thomas,
  eds.), Springer, Berlin, 2017, pp.~115--121.

\bibitem[BKP82]{Barany:1982ga}
I.~B{\'a}r{\'a}ny, M.~Katchalski, and J.~Pach, \emph{{Quantitative Helly-type
  theorems}}, Proc. Amer. Math. Society \textbf{86} (1982), 109--114.

\bibitem[BKP21]{BarKP:2020}
I.~B{\'a}r{\'a}ny, G.~Kalai, and A.~P\'or, \emph{{Universal sequences of lines
  in $\R^d$}}, 2021, arXiv:2110.12474, p.~21.

\bibitem[BL92]{Barany:1992tx}
I.~B{\'a}r{\'a}ny and D.~G. Larman, \emph{{A coloured version of Tverberg's
  theorem}}, J. London Math. Society \textbf{45} (1992), 314--320.

\bibitem[BLN17]{Bukh:2017vs}
B.~Bukh, Sh-P. Loh, and G.~Nivasch, \emph{{Classifying unavoidable Tverberg
  partitions}}, J. Comput. Geom. \textbf{8} (2017), 174--205.

\bibitem[BLVS{\etalchar{+}}93]{Bjor:1993}
A.~Bj\"orner, M.~Las~Vergnas, B.~Sturmfels, N.~White, and G.~M. Ziegler,
  \emph{Oriented matroids}, Encyclopedia Math. Appl., vol.~46, Cambridge
  University Press, Cambridge, 1993.

\bibitem[BM03]{Barany:2003wg}
I.~B{\'a}r{\'a}ny and J.~Matou{\v s}ek, \emph{{A fractional Helly theorem for
  convex lattice sets}}, Adv. Math. \textbf{174} (2003), 227--235.

\bibitem[BMN11]{Bukh:2011vs}
B.~Bukh, J.~Matou{\v s}ek, and G.~Nivasch, \emph{{Lower bounds for weak
  epsilon-nets and stair-convexity}}, Israel J. Math. \textbf{182} (2011),
  199--228.

\bibitem[BMNT18]{BarMesh:2018}
I.~B{\'a}r{\'a}ny, R.~Meshulam, E.~Nevo, and M.~Tancer, \emph{Pach's theorem
  does not admit a topological extension}, Discrete Comput. Geom. \textbf{360}
  (2018), 420--429.

\bibitem[BMP16]{BarMP:2016}
I.~B\'ar\'any, J.~Matou\v{s}ek, and A.~P{\'o}r, \emph{{Curves in ${\mathbb
  R}^d$ intersecting every hyperplane at most $d+1$ times}}, {J. European Math.
  Society} \textbf{11} (2016), 2469--2482.

\bibitem[BMZ15]{blago15}
P.~V.~M. Blagojevi{\'c}, B.~Matschke, and G.~M. Ziegler, \emph{Optimal bounds
  for the colored {T}verberg problem}, J. European Math. Society \textbf{17}
  (2015), 739--754.

\bibitem[BP90]{BarPer:90}
I.~B{\'a}r{\'a}ny and M.~Perles, \emph{{The Carath\'eodory number for the
  $k$-core}}, Combinatorica \textbf{10} (1990), 185--194.

\bibitem[Bra17]{Brazit:2017}
S.~Brazitikos, \emph{Quantitative {H}elly-type theorem for the diameter of
  convex sets}, Discrete Comput. Geometry \textbf{57} (2017), 494--505.

\bibitem[BSS81]{Barany:1981vh}
I.~B{\'a}r{\'a}ny, S.~B. Shlosman, and A.~Sz\H{u}cs, \emph{{On a topological
  generalization of a theorem of Tverberg}}, J. London Math. Society \textbf{2}
  (1981), 158--164.

\bibitem[Buk10]{bukh2010radon}
B.~Bukh, \emph{Radon partitions in convexity spaces}, 2010, arXiv:1009.2384,
  p.~11.

\bibitem[Car07]{Carath1907}
C.~Carath\'eodory, \emph{{\"Uber den Variabilit\"atsbereich der Koeffizienten
  von Potenzreihen}}, Math. Annalen \textbf{64} (1907), 95--115.

\bibitem[CFP{\etalchar{+}}14]{Conlon:2014fx}
D.~Conlon, J.~Fox, J.~Pach, B.~Sudakov, and A.~Suk, \emph{{Ramsey-type results
  for semi-algebraic relations}}, Trans. American Math. Society \textbf{366}
  (2014), 5043--5065.

\bibitem[Cho17]{Polymath12}
T.~Chow, \emph{Rota's basis conjecture: Polymath 12}, 2017,
  https://polymathprojects.org/2017/05/05/rotas-basis-conjecture-polymath-12-p%
ost-3/.

\bibitem[CSSS20]{ChSSS:2020}
M.~Chudnovsky, A.~Scott, P.~Seymour, and S.~Spirkl, \emph{Proof of the
  {K}alai--{M}eshulam conjecture}, Israel J. Math. \textbf{238} (2020),
  639--661.

\bibitem[Deb70]{Debrunner:1970th}
H.~Debrunner, \emph{{Helly type theorems derived from basic singular
  homology}}, Amer. Math. Monthly \textbf{77} (1970), 375--380.

\bibitem[DFN21]{DFN:2020}
G.~Dam\'asdi, V.~F\"oldv\'ari, and M.~Nasz\'odi, \emph{Colorful {H}elly-type
  theorems for the volume of intersections of convex bodies}, J. Combin. Theory
  A \textbf{178} (2021), 5043--5065.

\bibitem[DGHP05]{Dhand:2005}
R.~Dhandapani, J.~E. Goodman, A.~Holmsen, and R.~Pollack, \emph{Interval
  sequences and the combinatorial encoding of planar families of convex sets},
  Rev. Roumaine Math. Pures Appl. \textbf{50} (2005), 537--553.

\bibitem[DGK63]{Danzer:1963ug}
L.~Danzer, B.~Gr{\"u}nbaum, and V.~Klee, \emph{{Helly's theorem and its
  relatives}}, Proc. Sympos. Pure Math., Vol. VII, Amer. Math. Society,
  Providence, RI., 1963, pp.~101--180.

\bibitem[DLGMM19]{de2017discrete}
J.~A. De~Loera, X.~Goaoc, F.~Meunier, and N.~Mustafa, \emph{{The discrete yet
  ubiquitous theorems of Carath\'eodory, Helly, Sperner, Tucker, and
  Tverberg}}, Bull. American Math. Society \textbf{26} (2019), 415--511.

\bibitem[Doi73]{Doignon:1973ht}
J.~P. Doignon, \emph{{Convexity in cristallographical lattices}}, J. Geom.
  \textbf{3} (1973), 71--85.

\bibitem[DS21]{DillSob:2020}
T.~Dillon and P.~Sober\'on, \emph{{A m\'elange of diameter of Helly-type
  theorems}}, SIAM J. Discrete Math. \textbf{35} (2021), 1615--1627.

\bibitem[Eck79]{Eckhoff:1979bi}
J.~Eckhoff, \emph{Radon's theorem revisited}, Contributions to geometry
  ({P}roc. {G}eom. {S}ympos., {S}iegen, 1978), Birkh\"auser, Basel, 1979,
  pp.~164--185.

\bibitem[Eck85]{Eckhoff:1985hr}
\bysame, \emph{{An upper-bound theorem for families of convex sets}}, Geom.
  Dedicata \textbf{19} (1985), 217--227.

\bibitem[Eck88]{Eckhoff:1988vk}
\bysame, \emph{{Intersection properties of boxes. Part I: An upper-bound
  theorem}}, Israel J. Math. \textbf{62} (1988), 283--301.

\bibitem[Eck91]{Eckhoff:1991ue}
\bysame, \emph{{Intersection properties of boxes. Part II: Extremal families}},
  Israel J. Math. \textbf{73} (1991), 129--149.

\bibitem[Eck93]{Eck93survey}
\bysame, \emph{Helly, {R}adon, and {C}arath\'eodory type theorems}, Handbook of
  convex geometry, {V}ol.\ {A}, {B}, North-Holland, Amsterdam, 1993,
  pp.~389--448.

\bibitem[Eck00]{Eckhoff:2000jw}
\bysame, \emph{{The partition conjecture}}, Discrete Math. \textbf{221} (2000),
  61--78.

\bibitem[Eck03]{Eckhoff:2003ed}
\bysame, \emph{{A Survey of the Hadwiger--Debrunner (p,q)-problem}}, Discrete
  and Computational Geometry, Springer, Berlin, Berlin, Heidelberg, 2003,
  pp.~347--377.

\bibitem[EN09]{Eckhoff:2009kv}
J.~Eckhoff and K.-P. Nischke, \emph{{Morris's pigeonhole principle and the
  Helly theorem for unions of convex sets}}, Bull. Lond. Math. Society
  \textbf{41} (2009), 577--588.

\bibitem[FGL{\etalchar{+}}12]{FoxGr:2012}
J.~Fox, M.~Gromov, V.~Lafforgue, A.~Naor, and J.~Pach, \emph{{Overlap
  properties of geometric expanders}}, J. reine und angewandte Mathematik
  \textbf{671} (2012), 49--83.

\bibitem[FPT11]{Fox:2011ca}
J.~Fox, J.~Pach, and C.~D. T{\'o}th, \emph{{Intersection patterns of curves}},
  J. London Math. Society \textbf{83} (2011), 389--406.

\bibitem[Fri15]{Frick:2015}
F.~Frick, \emph{{Counterexamples to the topological Tverberg conjecture}},
  Oberwolfach Reports \textbf{12} (2015), 318--321.

\bibitem[FS20]{FriSob:2020}
F.~Frick and P.~Sober\'on, \emph{{The topological Tverberg beyond prime
  powers}}, 2020, arXiv:2005.05251, p.~13.

\bibitem[GLS08]{GaoLS:2008ej}
J.~Gao, M.~Langberg, and L.~J. Schulman, \emph{{Analysis of incomplete data and
  an intrinsic-dimension Helly theorem}}, Discrete Comput. Geom. \textbf{40}
  (2008), 537--560.

\bibitem[GM61]{Grunbaum:1961fd}
B.~Gr{\"u}nbaum and T.~S. Motzkin, \emph{{On components in some families of
  sets}}, Proc. American Math. Society \textbf{12} (1961), 607--613.

\bibitem[GP85]{Goodman:1993ji}
J.~E. Goodman and R.~Pollack, \emph{Allowable sequences and order types in
  discrete and computational geometry}, New Trends in Discrete and
  Computational Geometry (J.~Pach, ed.), Springer, Berlin, 1985, pp.~103--134.

\bibitem[GPP{\etalchar{+}}17]{GPPTW:2017dp}
X.~Goaoc, P.~Pat\'ak, Z.~Pat\'akova, M.~Tancer, and U.~Wagner, \emph{{Bounding
  Helly numbers via Betti numbers}}, Journey through Discrete Mathematics. A
  Tribute to Ji\v{r}\'{\i} Matou\v{s}ek (M.~Loebl, J.~Ne\v{s}et\v{r}il, and
  R.~Thomas, eds.), Springer, Berlin, 2017, pp.~407--447.

\bibitem[GPW93]{Goodman:1993xc}
J.~E. Goodman, R.~Pollack, and R.~Wenger, \emph{Geometric transversal theory},
  New Trends in Discrete and Computational Geometry, Springer, Berlin, 1993,
  pp.~163--198.

\bibitem[Gro10]{Gromov:2010eb}
M.~Gromov, \emph{{Singularities, expanders and topology of maps. Part 2: From
  combinatorics to topology via algebraic isoperimetry}}, Geometric Func. Anal.
  \textbf{20} (2010), 416--526.

\bibitem[HD57]{Hadwiger:1957we}
H.~Hadwiger and H.~Debrunner, \emph{{{\"U}ber eine Variante zum Hellyschen
  Satz}}, Arch. Math \textbf{8} (1957), 309--313.

\bibitem[Hel23]{Helly:1923wr}
E.~Helly, \emph{\"{U}ber {M}engen konvexer {K}\"orper mit gemeinschaftlichen
  {P}unkten}, Jahresberichte der Deutschen Math.-Verein. \textbf{32} (1923),
  175--176.

\bibitem[Hel30]{Helly:1930hk}
\bysame, \emph{{{\"U}ber Systeme von abgeschlossenen Mengen mit
  gemeinschaftlichen Punkten}}, Monatsh. f\"ur Mathematik und Physik
  \textbf{37} (1930), 281--302.

\bibitem[Hel07]{Hell:2007tp}
S.~Hell, \emph{{On the number of Tverberg partitions in the prime power case}},
  European J. Combin. \textbf{28} (2007), 347--355.

\bibitem[HL21]{Holm:2019uv}
A.~F. Holmsen and D-G. Lee, \emph{Radon numbers and the fractional {H}elly
  theorem}, Israel J. Math. \textbf{241} (2021), 433--447.

\bibitem[Hol13]{Holmsen:2013tu}
A.~F. Holmsen, \emph{{Geometric transversal theory: T(3)-families in the
  plane}}, Geometry: Intuitive, Discrete, and Convex, J\'anos Bolyai Math.
  Society, Budapest, 2013, pp.~187--203.

\bibitem[IN21]{NaIv:2021}
G.~Ivanov and M.~Nasz\'odi, \emph{Quantitative {H}elly-type theorem:
  containment in a homothet}, 2021, arXiv:2103.04122, p.~6.

\bibitem[Kal84a]{Kalai:1984isa}
G.~Kalai, \emph{{Characterization of $f$-vectors of families of convex sets in
  $R^d$ Part I: Necessity of Eckhoff{\textquoteright}s conditions}}, Israel J.
  Math. \textbf{48} (1984), 175--195.

\bibitem[Kal84b]{Kalai:1984bg}
\bysame, \emph{{Intersection patterns of convex sets}}, Israel J. Math.
  \textbf{48} (1984), 161--174.

\bibitem[Kal86]{Kalai:1986hoa}
\bysame, \emph{{Characterization of $f$-vectors of families of convex sets in
  $R^d$, part II: Sufficiency of Eckhoff's conditions}}, J. Combin. Theory,
  Ser. A \textbf{41} (1986), 167--188.

\bibitem[Kal95]{kalai:1994}
\bysame, \emph{{Combinatorics and convexity}}, Proceedings of the International
  Congress of Mathematicians, Z\"urich, Birkh\"auser, Basel, 1995,
  pp.~1363--1374.

\bibitem[Kal00]{Kal00flavor}
\bysame, \emph{Combinatorics with a geometric flavor}, Visions in Mathematics
  (N. Alon, J. Bourgain, M. Gromov, and V. Milman, eds.), Birkh\"auser, Basel,
  2000, pp.~742--791.

\bibitem[Kal02]{kalai:2002ij}
\bysame, \emph{Algebraic shifting, computational commutative algebra and
  combinatorics}, Adv. Stud. Pure Math., vol.~33, Math. Soc. Japan, Tokyo,
  2002, pp.~121--163.

\bibitem[Kal10]{kalai:2010}
\bysame, \emph{{ Combinatorial and topological aspects of Helly type
  theorems}}, https://gilkalai.files.wordpress.com/2010/10/es.pdf, 2010.

\bibitem[Kat71]{Katchalski:1971to}
M.~Katchalski, \emph{{The dimension of intersections of convex sets}}, Israel
  J. Math. \textbf{10} (1971), 465--470.

\bibitem[Kat78]{Katchalski:1978to}
\bysame, \emph{{Reconstructing dimensions of intersections of convex sets}},
  Aequationes Math. \textbf{17} (1978), 249--254.

\bibitem[KGT01]{Kleitman:2001vo}
D.~J. Kleitman, A.~Gy{\'a}rf{\'a}s, and G.~T{\'o}th, \emph{{Convex sets in the
  plane with three of every four meeting}}, Combinatorica \textbf{21} (2001),
  221--232.

\bibitem[Kim17]{Kim:2017}
M.~Kim, \emph{{Note on the colorful fractional Helly theorem}}, Discrete Math
  \textbf{340} (2017), 3167--3170.

\bibitem[KL79]{Katchalski:1979vt}
M.~Katchalski and A.~Liu, \emph{{A problem of geometry in ${\mathbb R}^n$}},
  Proc. Amer. Math. Society \textbf{75} (1979), 284--288.

\bibitem[KM05]{Kalai:2005cm}
G.~Kalai and R.~Meshulam, \emph{{A topological colorful Helly theorem}}, J.
  Combin. Theory, Ser. A \textbf{191} (2005), 305--311.

\bibitem[KM08]{Kalai:2008kc}
\bysame, \emph{{Leray numbers of projections and a topological Helly-type
  theorem}}, J. Topology \textbf{1} (2008), 551--556.

\bibitem[KS18]{KellSmor:2018}
Ch. Keller and Sh. Smorodinsky, \emph{{On piercing number of families
  satisfying the $(p,q)_r$ property}}, Comput. Geometry \textbf{72} (2018),
  11--18.

\bibitem[KS20]{KellSmor:2020}
\bysame, \emph{{From a (p, 2)-theorem to a tight (p, q)-theorem}}, Discrete
  Comput. Geometry \textbf{63} (2020), 821--847.

\bibitem[KST18]{KellSmoTar:2018}
Ch. Keller, Sh. Smorodinsky, and G.~Tardos, \emph{{Improved bounds on the
  Hadwiger--Debrunner numbers}}, Israel J. Math. \textbf{225} (2018), 925--945.

\bibitem[Lar68]{Larman:1968}
D.~G. Larman, \emph{{Helly type properties of unions of convex sets}},
  Mathematika \textbf{15} (1968), 53--59.

\bibitem[Lar72]{Larman:1072}
\bysame, \emph{On sets projectively equivalent to the vertices of a convex
  polytope}, Bull. London Math. Soc. \textbf{4} (1972), 6--12.

\bibitem[LB62]{LeBo:1962}
C.~G. Lekkerkerker and J.~Ch. Boland, \emph{{Representation of a finite graph
  by a set of intervals on the real line}}, Fund. Math. \textbf{51} (1962),
  45--64.

\bibitem[LMPT94]{Larman:1994wt}
D.~G. Larman, J.~Matou{\v s}ek, J.~Pach, and J.~T\"or\H{o}csik, \emph{{A
  Ramsey-type result for convex sets}}, Bull. London Math. Society \textbf{26}
  (1994), 132--136.

\bibitem[Mat97]{Matousek:1997di}
J.~Matou{\v s}ek, \emph{{A Helly-type theorem for unions of convex sets}},
  Discrete Comput. Geom. \textbf{18} (1997), 1--12.

\bibitem[Mat04]{Matousek:2004cs}
\bysame, \emph{{Bounded VC-dimension implies a fractional Helly theorem}},
  Discrete Comput. Geom. \textbf{31} (2004), 251--255.

\bibitem[McG20]{McGinn:2020}
D.~McGinnis, \emph{A family of convex sets in the plane satisfying the
  $(4,3)$-property can be pierced by nine points}, 2020, arXiv:2010.13195,
  p.~16.

\bibitem[McM70]{McM70}
P.~McMullen, \emph{{The maximum numbers of faces of a convex polytope}},
  Mathematika \textbf{17} (1970), 179--184.

\bibitem[MMSS17]{MeunM:2017}
F.~Meunier, W.~Mulzer, P.~Sarabezolles, and Y.~Stein, \emph{{The rainbow at the
  end of the line: A PPAD formulation of the colorful Carath\'eodory theorem
  with applications}}, Proceedings of the Twenty-Eighth {A}nnual {ACM}-{SIAM}
  {S}ymposium on {D}iscrete {A}lgorithms, ACM, New York, 2017, pp.~1342--1351.

\bibitem[Mon14]{Montejano:2014ii}
L.~Montejano, \emph{{A new topological Helly theorem and some transversal
  results}}, Discrete Comput. Geom. \textbf{52} (2014), 390--398.

\bibitem[Mot55]{Motzkin:1955}
T.~S. Motzkin, \emph{{A proof of Hilbert's Nullstellensatz}}, Math. Zeit.
  \textbf{63} (1955), 341--344.

\bibitem[MRR20]{MarRR:2018}
L.~Mart\'inez, E.~Rold\'an, and N.~Rubin, \emph{{Further consequences of the
  colorful Helly hypothesis.}}, Discrete Computat. Geometry \textbf{63} (2020),
  848--866.

\bibitem[MS11]{Montejano:2011cz}
L.~Montejano and P.~Sober\'on, \emph{{Piercing numbers for balanced and
  unbalanced families}}, Discrete Comput. Geom. \textbf{45} (2011), 358--364.

\bibitem[MW14]{MW14SoCG}
I.~Mabillard and U.~Wagner, \emph{Eliminating {T}verberg points, {I}. {A}n
  analogue of the {W}hitney trick}, Computational Geometry ({S}o{CG}'14), ACM,
  New York, 2014, pp.~171--180.

\bibitem[MY20]{MorYe:2020}
S.~Moran and A.~Yehudayoff, \emph{On weak $\eps$-nets and the {R}adon number},
  Discrete Comput. Geom. \textbf{64} (2020), 1125--1140.

\bibitem[Nas16]{Naszodi:2015vi}
M.~Nasz{\'o}di, \emph{{Proof of a conjecture of B{\'a}r{\'a}ny, Katchalski and
  Pach}}, Discrete Comput. Geom. \textbf{55} (2016), 243--248.

\bibitem[Onn01]{Onn:2001}
Sh. Onn, \emph{The {R}adon-split and the {H}elly-core of a point
  configuration}, J. of Geometry \textbf{72} (2001), 157--162.

\bibitem[{\"O}za87]{Oza87}
M.~{\"O}zaydin, \emph{{Equivariant maps for the symmetric group}}, unpublished
  preprint, University of Winsconsin-Madison, 17 pages, 1987.

\bibitem[Pac98]{Pach:1998vx}
J.~Pach, \emph{{A Tverberg-type result on multicolored simplices}}, Comput.
  Geom. \textbf{10} (1998), 71--76.

\bibitem[P{\'a}l20]{Palv2020}
D.~P{\'a}lv\"olgyi, \emph{Radon numbers grow linearly}, 36th international
  Symposium on Computational Geometry, vol. 164, Schloss Dagstuhl.
  Leibniz-Zent. Inform., Wadern, 2020, LIPIcs. Leibniz Int. Proc. Inf., p.~5.

\bibitem[P{\'o}r18]{pornew}
A.~P{\'o}r, \emph{{Universality of vector sequences and universality of
  Tverberg partitions}}, 2018, arXiv:1805.07197, p.~30.

\bibitem[PS16]{PeS16}
M.~A. Perles and M.~Sigron, \emph{Some variations on {T}verberg's theorem},
  Israel J. Math. \textbf{216} (2016), 957--972. \MR{3557472}

\bibitem[RA01]{Ram01}
J.~L. Ram\'irez-Alfons\'in, \emph{Lawrence oriented matroids and a problem of
  {M}c{M}ullen on projective equivalences of polytopes}, European J. Combin.
  \textbf{22} (2001), 723--731.

\bibitem[Rad21]{Radon:1921vh}
J.~Radon, \emph{{Mengen konvexer K{\"o}rper, die einen gemeinsamen Punkt
  enthalten}}, Math. Ann. \textbf{83} (1921), 113--115.

\bibitem[Rea79]{Reay79}
J.~R. Reay, \emph{Several generalizations of {T}verberg's theorem}, Israel J.
  Math. \textbf{34} (1979), 238--244.

\bibitem[Rub18]{Ru18}
N.~Rubin, \emph{An improved bound for weak $\varepsilon$-nets in the plane},
  Proceedings of the Annual Symposium on Foundations of Computer Science
  ({FOCS}), 2018, pp.~224--235.

\bibitem[Rub21]{Ru20}
\bysame, \emph{Stronger bounds for weak $\varepsilon$-nets in higher
  dimensions}, Proceedings of the Annual Symposium on Foundations of Computer
  Science ({STOC 2021}), 2021, arxiv:2104.12654, p.~62.

\bibitem[Sar92]{Sarkaria:1992vt}
K.~S. Sarkaria, \emph{{Tverberg{\textquoteright}s theorem via number fields}},
  Israel J. Math. \textbf{79} (1992), 317--320.

\bibitem[Sie79]{Sierksma:1979}
G.~Sierksma, \emph{Convexity without linearity; the {D}utch cheese problem},
  Mimeographed notes, 1979.

\bibitem[Sob15]{Sober:2015}
P.~Sober\'on, \emph{{Equal coefficients and tolerance in coloured Tverberg
  partitions}}, Combinatorica \textbf{35} (2015), 235--252.

\bibitem[SS20]{ScottSey:2020}
A.~Scott and P.~Seymor, \emph{A survey of $\chi$-boundedness}, J. Graph Theory
  \textbf{95} (2020), 473--504.

\bibitem[Sta75]{Stanley:1975}
R.~P. Stanley, \emph{The upper bound conjecture and {C}ohen—-{M}acaulay
  rings}, Studies in Appl. Math. \textbf{54} (1975), 135--142.

\bibitem[Suk14]{Suk:2014tt}
A.~Suk, \emph{{A note on order-type homogeneous point sets}}, Mathematika
  \textbf{60} (2014), 37--42.

\bibitem[Tan13]{Tancer:2013iz}
M.~Tancer, \emph{{Intersection patterns of convex sets via simplicial
  complexes: A survey}}, Thirty essays on geometric graph theory, Springer, New
  York, NY, 2013, pp.~521--540.

\bibitem[Tve66]{Tverberg:1966tb}
H.~Tverberg, \emph{{A generalization of Radon{\textquoteright}s theorem}}, J.
  London Math. Society \textbf{41} (1966), 123--128.

\bibitem[V{\v Z}92]{ZivV:1992}
S.~T. Vre{\'c}ica and R.~T. {\v Z}ivaljevi{\'c}, \emph{{The colored Tverberg's
  problem and complexes of injective functions}}, J. Combin. Theory Ser. A
  \textbf{61} (1992), 309--318.

\bibitem[V{\v Z}93]{Vucic:1993be}
A.~Vu{\v c}i{\'c} and R.~T. {\v Z}ivaljevi{\'c}, \emph{{Note on a conjecture of
  {S}ierksma}}, Discrete Comput. Geom. \textbf{9} (1993), 339--349.

\bibitem[Weg75]{Wegner:1975eo}
G.~Wegner, \emph{{d-collapsing and nerves of families of convex sets}}, Arch.
  Math. \textbf{26} (1975), 317--321.

\bibitem[Wen99]{Wenger:2001}
R.~Wenger, \emph{Progress in geometric transversal theory}, Advances in
  discrete and computational geometry, Contemp. Math., vol. 223, American Math.
  Society, Providence, RI, 1999, pp.~375--393.

\bibitem[Whi17]{Whi17}
M.~J. White, \emph{On {T}verberg partitions}, Israel J. Math. \textbf{219}
  (2017), 549--553.

\bibitem[Whi21]{White2020}
\bysame, \emph{{A new topological property of nerves of convex sets in
  ${\mathbb R}^d$}}, manuscript, 2021.

\end{thebibliography}

\bigskip

\noindent
Imre B\'ar\'any \\
R\'enyi Institute of Mathematics,\\
13-15 Re\'altanoda Street, Budapest, 1053 Hungary\\
{\tt barany.imre@renyi.hu} and\\
Department of Mathematics\\
University College London\\
Gower Street, London, WC1E 6BT, UK

\medskip
\noindent
Gil Kalai\\
Einstein Institute of Mathematics\\
Hebrew University,
Jerusalem 91904, Israel,\\
{\tt kalai@math.huji.ac.il} and\\
Efi Arazy School of Computer Science,
IDC, Herzliya, Israel

\end{document}